\documentclass{elsarticle}

\usepackage{lineno,hyperref}
\usepackage{subcaption}
\usepackage[justification=centering]{caption}
\usepackage{graphicx}
\graphicspath{ {./figures/} }
\usepackage{siunitx}
\usepackage{numprint}
\usepackage[dvipsnames]{xcolor}
\usepackage[normalem]{ulem}
\usepackage{amsmath}
\usepackage{amsthm}
\usepackage{amssymb}
\usepackage{blindtext}
\usepackage{nicefrac}
\usepackage{multicol}
\usepackage{bm}
\usepackage{blindtext}
\usepackage{cancel} % this may not be used. Only used for cancelling arrow.

\usepackage[utf8]{inputenc}
\modulolinenumbers[5]

\definecolor{darkgreen}{rgb}{0.0, 0.35, 0.0}
\newcommand{\reva}[1]{\textcolor{black}{#1}}
\newcommand{\revb}[1]{\textcolor{black}{#1}}

\newcommand{\bfdot}[1]{\dot{\bm{#1}}}

\newcommand{\DDt}[1]{\frac{D\bm{#1}}{Dt}}

\journal{}

\begin{document}

\begin{frontmatter}

\title{Efficient numerical methods for the Maxey-Riley-\reva{Gatignol} equations with Basset history term\tnoteref{mytitlenote}}
\tnotetext[mytitlenote]{This project is funded by the Deutsche Forschungsgemeinschaft (DFG, German Research Foundation) – SFB 1615 – 503850735.}

%% Group authors per affiliation:
\author[tuhh]{Julio Urizarna-Carasa}
\ead{julio.urizarna@tuhh.de}

\author[tuhh]{Leon Schlegel}
\ead{leon.schlegel@tuhh.de}

\author[tuhh]{Daniel Ruprecht\corref{cor}}
\ead{ruprecht@tuhh.de}

\address[tuhh]{Lehrstuhl Computational Mathematics, Institut für Mathematik, Technische Universität Hamburg, Germany}

\cortext[cor]{Corresponding author}

\begin{abstract}
The Maxey-Riley\reva{-Gatignol} equations (MRGE) describe the motion of a finite-sized, spherical particle in a  fluid.
Because of wake effects, the force acting on a particle depends on its past trajectory.
This is modelled by an integral term in the MRGE, also called Basset force, that makes its numerical solution challenging and memory intensive.
A recent approach proposed by \revb{[S. G. Prasath, V. Vasan, R. Govindarajan, Accurate solution method for the Maxey–Riley equation, and the effects of Basset history, Journal of
Fluid Mechanics 868 (2019) 428–460]} exploits connections between the integral term and fractional derivatives to reformulate the MRGE as a time-dependent partial differential equation on a semi-infinite pseudo-space.
They also propose a numerical algorithm based on polynomial expansions.
This paper develops a numerical approach based on finite difference instead, by adopting techniques by \revb{[M. N. Koleva, Numerical solution of the heat equation in unbounded do-
mains using quasi-uniform grids, in: International Conference on Large-
Scale Scientific Computing, Springer, 2005, pp. 509–517]} \revb{and} \revb{[R. Fazio, A. Jannelli, Finite difference schemes on quasi-uniform grids for BVPs on infinite intervals, Journal of Computational and Applied Mathematics 269 (2014) 14–23]} to cope with the issues of having an unbounded spatial domain.
We compare convergence order and computational efficiency for particles of varying size and density of the polynomial expansion by Prasath et al., our finite difference schemes and a direct integrator for the MRGE based on multi-step methods proposed by \revb{[A. Daitche, Advection of inertial particles in the presence of the history force: Higher order numerical schemes, Journal of Computational Physics 254 (2013) 93–106]}.
\reva{While all methods achieve their theoretical convergence order for neutrally buoyant particles with zero initial relative velocity, they suffer from various degrees of order reduction if the initial relative velocity is non-zero or the particle has a different density than the fluid.}
\end{abstract}


\end{frontmatter}

\section{Introduction}
The Maxey-Riley\reva{-Gatignol} equations (MRGE) are a second-order system of integro-differential equations with a singular kernel that model the motion of finite-sized (inertial), spherical particles in a fluid~\cite{MaxeyandRiley1983}.
They are used for the study of a wide range of phenomena, for example the formation of clouds in the atmosphere~\cite{GovindarajanAndRavichandran2017, FalkovichEtAl2002}, the settling of plankton in the ocean (``marine snow'')~\cite{NiaziEtAl2016, GusevaEtAl2016}, the transmission of Covid-19 virus through sprays~\cite{CumminsEtAl2020} and more.

While the motion of inertial particles can be described by other means, for example the Navier-Stokes equations with specific boundary conditions~\cite{GaldiEtAl2008, CartwrightEtAl2010}, the MRGE provide a  simpler and computationally much more efficient model which is valid for spheric particles under the assumption that while the flow influences the particle, the particles do not alter the fluid. 
This is an appropriate assumption when ``particles are far smaller than all relevant flow length scales, and when the suspension of particles is dilute''~\cite{PrasathEtAl2019}. 
The MRGE are useful for particles that are too large to be treated as passive tracers but not yet big enough to substantially disturb the fluid. 
Candelier et al.~\cite[Figure 3]{CandelierEtAl2004} show a good match between theoretical and experimental trajectories for particles in a rotating cylinder.

According to Langlois et al.~\cite{LangloisEtAl2015}, the search for a suitable model that ultimately resulted in the MRGE dates back to 1851, when Stokes attempted for the first time to obtain a model for the motion of a pendulum in a fluid~\cite{Stokes1851}. 
Later, the combined works of Basset~\cite{Basset1888}, Boussinesq~\cite{Boussinesq1885} and Oseen~\cite{Oseen1927} let to a model called the Basset-Boussinesq-Oseen (BBO) equation for the settling of a sphere under gravity in a quiescent fluid.
In 1947, Tchen~\cite{Tchen2013} modified the BBO into a model for the motion of a rigid sphere in a nonuniform unsteady flow. 
This model underwent several amendments until it attained the form that is now known as the Maxey-Riley\reva{-Gatignol} equations~\cite{MaxeyandRiley1983,Gatignol1983}.
Throughout this paper, we consider the MRGE in the form 
\begin{subequations}
    \label{eq:MR}
\begin{align}
		\bfdot{y}(t) &= \bm{v}(t), \label{MR1} \\
		R\bfdot{v}(t) &= \DDt{u} - \frac{1}{S} (\bm{v} - \bm{u}) - \sqrt{\frac{3}{\pi S}} \left\{ \frac{1}{\sqrt{t}} \left(\bm{v}(0) - \bm{u}(0)\right) + \int_0^t \frac{(\bfdot{v}(s) - \bfdot{u}(s))}{\sqrt{t-s}}ds \right\}, \label{MR2}
\end{align}
\end{subequations}
given by Prasath et al.~\cite{PrasathEtAl2019} which neglects the Faxen corrections by assuming that the particle radius is much smaller than the typical length scale of the flow~\cite{BeronEtAl2019}.
Here, $\bm{u}(\bm{y}(t),t)$ is the fluid velocity at the particle position $\bm{y}(t)$ and $\bm{v}(t)$ is the particle velocity. 
The parameters are
\begin{align}
    \label{eq:parameters}
	\beta = \frac{\rho_p}{\rho_f}, \quad R = \frac{1+2\beta}{3}, \quad S = \frac{a^2 }{3T\nu},
\end{align}
where $\rho_p$ is the density of the particle, $\rho_f$ the density of the fluid, $a$ the particle's diameter, $\nu$ the fluid's kinematic viscosity and $T$ the time scale associated with the flow.

Due to the integral term in~\eqref{MR2}, also known as Basset history or memory term, the force acting on a particle depends on its past trajectory and the MRGE are not a dynamical system.
Since this term causes numerical difficulties and, in straightforward numerical integration, requires storage of all previously computed steps, it is often neglected in applications~\cite{CumminsEtAl2020}, modified~\cite{LovalentiEtAl1993, Mei1994, DorganAndLoth2007, MorenoCasasEtAl2016} or approximated~\cite{KlinkenbergEtAl2014, ElghannayAndTafti2016, ParmarEtAl2018, VanHinsbergEtAl2011}. 
However, both theoretical~\cite{PrasathEtAl2019} and empirical~\cite{CandelierEtAl2004, NinoAndGarcia1998, MordantAndPinton2000} studies have shown that the Basset history term is relevant even for small Stokes numbers and that neglecting it can lead to noticeable inaccuracies in the modelled trajectories.

Therefore, a handful of numerical approaches were developed to solve the full MRGE efficiently.
Based on linear multistep methods, Daitche proposed the first direct integration methods that do not rely on approximations of the kernels~\cite{Daitche2013}.
While computationally fast, the approach still requires storage of all time steps which can lead to problems with memory in long simulations.
Michaelides~\cite{michaelides1992} propose a transformed version of the MRGE based on the Laplace transform that removes the implicit dependency on the dependent variable in the history integral. 
While this avoids the need to use a nonlinear iterative solver and  thus reduces computing time, it does not remove the history integral.

Tatom in 1988 realized that the history term is equivalent to a half derivative of Riemann-Liouville type~\cite{Tatom1988}. 
Using this, Prasath et al.~\cite{PrasathEtAl2019} show that solving the MRGE is equivalent to solving a diffusive heat equation posed on a semi-infinite pseudo-space with a nonlinear boundary condition.
Their reformulation removes the integral term and memory effect and allows to use standard numerical techniques for partial differential equations to solve the MRGE.
However, it requires dealing with an unbounded spatial computational domain.
Prasath et al.~\cite{PrasathEtAl2019} also propose a numerical approach based on the integral form of the solution obtained by using Fokas's method~\cite{Fokas1997}.
An open-source reimplementation in Python of their numerical approach based on polynomial expansions was recently published~\cite{UrizarnaEtAl2023}.
However, their reformulation also opens up possibilities to apply other numerical techniques for partial differential equations to the MRGE.
Very recently, Jaganathan et al.~\cite{JaganathanEtAl2023} proposed another transformation of the Maxey-Riley-Gatignol equation by embedding it into a an extended state-space.
This removes the non-locality in time and results in a dynamical system that can be solved with standard explicit numerical integrators.
\reva{While we did not reimplement their Fortran code in Python, a preliminary comparison of performance based on published results can be found in \S\ref{subsec:efficiency}. Comparing their method against Prasath et al's and our approach would be an interesting next step.}

To illustrate the opportunities offered by Prasath et al.'s~\reva{\cite{PrasathEtAl2019}} reformulation, we propose a finite difference method as solver for the MRGE.
By modifying techniques developed by Koleva~\cite{Koleva2005}, Alshina et al.~\cite{alshina2002numerical} and Fazio and Janelli~\cite{FazioAndJanelli2014} for solving initial-boundary value problems on infinite domains, we introduce a second and fourth order finite difference discretization.
To efficiently deal with the nonlinearity at the boundary, we use implicit-explicit Runge-Kutta methods of order two and four, which avoid the need for an iterative nonlinear solver in every stage.
We also provide an open-source Python implementation~\cite{juliouri_2024_10839777} of our finite difference approach, Prasath et al.'s algorithm~\reva{\cite{PrasathEtAl2019}} and Daitche's method~\reva{\cite{Daitche2013}}.
The code can be used to reproduce our results or to explore regimes and parameter configurations beyond what is reported in the paper.
To our knowledge, these are the only available algorithms that solve the MRGE without approximations to the kernel.
We perform a comprehensive comparison of all three methods with respect to accuracy and computational efficiency for five different flow fields and particles of different size and density.

\reva{We demonstrate that all three methods reach their full theoretical order of convergence only if the particle is neutrally buoyant and initially has zero relative velocity.
For particles that are not neutrally buoyant, we observe order reduction for all methods.
Even stronger order reduction is observed if the particle has non-zero initial relative velocity.
It is noteworthy that all methods experience this, albeit to varying degrees, even though they rely on different discretizations of different formulations of the MRGE.}

\revb{\paragraph{Contributions}
The main novel contributions of this paper are (i) the application of techniques for finite differences on semi-infinite domains to the reformulated MRGE, (ii) the idea to use an implicit-explicit Runge-Kutta method to avoid having to use a nonlinear solver, (iii) a comparison of different numerical methods for the full MRGE with respect to computational efficiency and (iv) a demonstration that the numerical solution of the full MRGE including history term in real-time is possible.
}

\revb{\paragraph{Structure of the paper}
Section~\ref{section:mre} summarizes the transformation of the MRGE into a diffusion-type PDE on a semi-infinite domain and introduces the finite differences used to solve the resulting system.
Section~\ref{section:results} first introduces the five used benchmark problems.
It then studies convergence rates of the investigated methods for particles of different densities.
Finally the section compares different methods with respect to computational efficiency by investigating the wallclock times required to reach a certain error.
Section~\ref{section:conclusions} summarizes the paper and draws conclusions.
}

\section{Numerical solution of the reformulated Maxey-Riley equations}\label{section:mre}
Subsection~\eqref{subsec:transform} briefly summarizes the reformulation of the MRGE with history term as a diffusive PDE on a semi-infinite domain~\cite{PrasathEtAl2019}.
Subsection~\eqref{subsec:fd} introduces a second-order and fourth-order finite difference scheme to discretize this PDE in a method-of-lines approach.
Subsection~\ref{subsection:time-stepping} describes two types of Runge-Kutta time stepping methods, diagonally implicit (DIRK) and implicit-explicit (IMEX), to efficiently integrate the resulting semi-discrete problem in time.

\subsection{Transforming the MRGE with Basset history}\label{subsec:transform}
We simply state the resulting transformed system here.
For a detailed explanation of the reformulation, we refer to the original paper by Prasath et al.~\reva{\cite{PrasathEtAl2019}} and our description of a reimplementation of their approach in an open-source Python software~\cite{UrizarnaEtAl2023}.
The reformulated problem reads
\begin{subequations}
\label{eq:MRSys}
\begin{align}
		\bm{q}_t(x,t) &= \bm{q}_{xx}(x,t), & x>0, t\in(0,T], \label{MRSys1} \\
		\bm{q}(x,0) &= \bm{0}, & x>0, \label{MRSys2}\\
		\bm{q}_t(0,t) + \alpha\bm{q}(0,t) - \gamma\bm{q}_x(0,t) &= \bm{f}(\bm{q}(0,t),\bm{y}(t),t), & t\in[0,T], \label{MRSys3} \\
		\bfdot{y}(t) &= \bm{q}(0,t) + \bm{u}(\bm{y}(t),t), & t\in[0,T], \label{MRSys4} \\
		\lim_{t\to0}\bm{q}(0,t) &= \bm{v}_0 - \bm{u}_0, \label{MRSys5}\\
		\bm{y}(0) &= \bm{y}_0. \label{MRSys6}
\end{align}
\end{subequations}
where $x$ is a pseudo-space without physical meaning, $t$ is time, $\bm{y}(t)$ is, as in~\eqref{eq:MR}, the particle's position at time $t$, and $\bm{q}(x,t)$ is a function with identical dimensionality as $\bm{y}(t)$, that is two for a two-dimensional flow field and three for a 3D field.
The boundary value $\bm{q}(0,t)$ is equal to the relative velocity of the particle at time $t$ and $\dot{\bm{y}}(t)$ is its absolute velocity. 
Furthermore, we have the right hand side function
\begin{equation}
    \label{eq:f}
    \bm{f}\left( \bm{q}(0,t),\bm{y}(t),t \right) := \left( \frac{1}{R} - 1 \right) \DDt{u} - \bm{q}(0,t) \cdot \nabla_y \bm{u}(\bm{y}(t),t)
\end{equation}
for the boundary condition and physical parameters
\begin{align}
    \label{eq:alpha&gamma}
    \alpha := \frac{1}{RS}, \quad \gamma := \frac{1}{R}\sqrt{\frac{3}{S}}
\end{align}
with $R$ and $S$ defined in~\eqref{eq:parameters}.
Solving~\eqref{eq:MRSys} produces the solution $\bm{y}(t)$ of~\eqref{eq:MR} with the advantage that the history term is not present.

\subsection{Finite differences for the transformed MRGE}\label{subsec:fd}
The transformed problem~\eqref{eq:MRSys} is defined on a semi-infinite computational domain, which requires some care when dealing with the \reva{boundary condition at infinity}.
To discretize equation~\eqref{eq:MRSys} in space, we propose two finite difference schemes: the second order scheme by Koleva~\cite{Koleva2005} based on the work by Alshina et al.~\cite{alshina2002numerical}, and a novel fourth order approximation obtained by using compact finite differences for uniform grids~\cite{lele1992} plus a new technique based on the core idea in Fazio and Jannelli~\cite{FazioAndJanelli2014} to map equidistant nodes to the semi-infinite domain.

We discretize the spatial domain with a quasi-uniform grid by defining a set of $N$ uniform grid points $\xi_n := \frac{n}{N}$, $n\in \{0,1,\dots,N-1 \}$ in the interval $[0,1)$.
These are then mapped to $[0, \infty)$ via the logarithmic mapping
\begin{align} \label{LogRule}
	x_n := x(\xi_n) = -c\ln(1-\xi_n),
\end{align}
where $c$ is a parameter that controls the distribution of nodes such that approximatelly half of the grid points are placed within the interval $[0,c]$.
In contrast to the algebraic rule also proposed by Koleva~\reva{\cite{Koleva2005}} and Fazio~\reva{\cite{FazioAndJanelli2014}}, the logarithmic rule produces a higher density of nodes around $x_0=0$.
Since our aim is to approximate the boundary value $\bm{q}(0,t)$, this is where we require the highest accuracy.
\reva{We have, however, not performed a detailed comparison of the two mappings.
Furthermore, while numerical experiments not documented here suggest that the value of $c$ can have a substantial impact on the overall error, we could not establish a robust heuristic how to choose it.
Therefore, we set $c=20$ in all numerical experiments, which provided satisfactory results, and abstained from tuning.
Significant further gains in accuracy might be possible using a more informed and may be even adaptive choice where $c$ changes over time but studying this is left for future work.}

Discretizing the spatial derivatives in~\eqref{MRSys1}-\eqref{MRSys3} together with~\eqref{MRSys4} results in the semidisrete system
\begin{align} \label{semi-discrete:sys}
    \underbrace{
    \begin{bmatrix}
        \dot{\bm{q}}(t) \\
        \dot{\bm{y}}(t)   
    \end{bmatrix}}_{=:\dot{\bm{\eta}}(t)}
    = \underbrace{
    \left[
    \begin{array}{c | c}
            A_s & 
            \begin{array}{c c}
            0 & 0 \\
            \vdots & \vdots \\
            0 & 0 \\
            \end{array} \\         
            \hline \\
            \begin{array}{c c c c c}
            1 & 0 & 0 & \dots & 0\\
            0 & 1 & 0 & \dots & 0
            \end{array}
            & 
            \begin{array}{c c}
            0 & 0 \\
            0 & 0
            \end{array}
    \end{array}
    \right]}_{=:A}
    \underbrace{
    \begin{bmatrix}
        \bm{q}(t) \\
        \bm{y}(t)   
    \end{bmatrix}}_{=:\bm{\eta}(t)}
    +
    \underbrace{
    \begin{bmatrix}
        \bm{v}(\bm{q}_0(t),\bm{y}(t),t) \\
        \bm{u}(\bm{y}(t),t)
    \end{bmatrix}}_{=:\bm{\omega}(\bm{q}_0(t),\bm{y}(t),t)}
\end{align}

where the specific forms of $A_s\in\mathbb{R}^{N\times N}$ and $\bm{v}(\bm{q}_0(t),\bm{y}(t),t)\in\mathbb{R}^N$ depend on whether the second or fourth order discretization described below is used.
Furthermore,
\begin{align} \label{q:vector}
    \bm{q}(t) :=
        \begin{bmatrix}
        q_0^{(1)}(t) & q_0^{(2)}(t) & q_1^{(1)}(t) & q_1^{(2)}(t) & \dots & q_{N-2}^{(1)}(t) & q_{N-2}^{(2)}(t)
        \end{bmatrix}^T,
\end{align}
 is a vector with both the horizontal and vertical components of the relative velocity at each node and
 \begin{align} \label{q0:vector}
    \bm{q}_0(t) :=
        \begin{bmatrix}
        q_0^{(1)}(t) & q_0^{(2)}(t)
        \end{bmatrix}^T.
\end{align}
Note that in a three-dimensional flow field we would have components $q^{(3)}_0(t)$, $q^{(3)}_1(t)$, etc. 
However, since we provide only numerical examples for 2D flow fields, we just show the two component case.
Below, we describe a second-order and fourth-order accurate approach to construct $A_s$.

\subsubsection{Second order spatial discretization}
\label{Second:order}
Koleva~\reva{\cite{Koleva2005}} provides the following derivative approximations for problems on one-dimenstional semi-infinite spatial domains
\begin{gather}
	\frac{\partial q(\xi(x),t)}{\partial x}\biggr\rvert_{x_{n+\nicefrac{1}{2}}} \approx \; \frac{q_{n+1}(t) - q_{n}(t)}{2\zeta_n}, \label{FD2:1stDeriv} \\
	\frac{\partial^2q(\xi(x),t)}{\partial x^2}\biggr\rvert_{x_n} \approx \; \frac{1}{\psi_n} \left[ \frac{\partial q(\xi(x),t)}{\partial x} \biggr\rvert_{x_{n+\nicefrac{1}{2}}} - \frac{\partial q(\xi(x),t)}{\partial x} \biggr\rvert_{n-\nicefrac{1}{2}} \right] \label{FD2:2ndDeriv}
\end{gather}
where $q$ corresponds to one of the components of the vector $\bm{q}$ used in~\eqref{eq:MRSys} and
\begin{align*}
\psi_n := x_{n+\nicefrac{1}{2}}-x_{n-\nicefrac{1}{2}}, \quad \zeta_n:=x_{n+\nicefrac{3}{4}} - x_{n+\nicefrac{1}{4}}
\end{align*}
for $n = 1, 2, ..., N-1$ are the mesh spacings in the pseudo-space coordinate.
At the left boundary we use
\begin{align} \label{FD2:Boundary}
    \frac{\partial q(\xi(x),t)}{\partial x}\biggr\rvert_{x_0} \approx \frac{1}{2} \left[ \frac{\partial q(\xi(x),t)}{\partial x}\biggr\rvert_{x_{\nicefrac{1}{2}}} + \frac{\partial q(\xi(x),t)}{\partial x}\biggr\rvert_{x_{-\nicefrac{1}{2}}} \right].
\end{align}
The boundary condition~\eqref{MRSys3} is coupled to the PDE~\eqref{MRSys1} by using the ghost grid point at $-\nicefrac{1}{2}$ in equations~\eqref{FD2:Boundary} and~\eqref{FD2:2ndDeriv} for $n=0$.
With these approximations, the entries $(a_{i,j})_{1\leq i \leq 2N,\; 1\leq j \leq 2N}$ of $A_s$ are
\begin{align}
    a_{11} =& \; a_{22} = - \frac{\gamma + 2\alpha\zeta_0}{\zeta_0(2+\gamma\psi_0)}, & a_{12} =& \; a_{24} = \frac{\gamma}{\zeta_0(2+\gamma \psi_0)}
\end{align}
for the boundary values and
\begin{align}
    a_{i,i-2} =& \; a_{i+1,i-1} = \frac{1}{2\psi_{\frac{i-1}{2}}\zeta_{\frac{i-3}{2}}}, \\
    a_{i,i} =& \; a_{i+1,i+1} = -\frac{1}{\psi_{\frac{i-1}{2}}} \left( \frac{\zeta_{\frac{i-3}{2}} + \zeta_{\frac{i-1}{2}}}{2\zeta_{\frac{i-3}{2}}\zeta_{\frac{i-1}{2}}}\right), \\
    a_{i,i+2} =& \; a_{i+1,i+3} = \frac{1}{2\psi_{\frac{i-1}{2}}\zeta_{\frac{i-1}{2}}},
\end{align}
for $i = 2m + 1$ with $m\in\{1,2,\dots, N-2\}$. 
All other matrix entries are zero so that $A_s$ is a sparse pentadiagonal matrix with zero first upper and lower diagonals. 
The vector $\bm{v}(\bm{q}_0(t),\bm{y}(t),t)$, which contains the boundary condition function $\bm{f}(\bm{q}_0(t),\bm{y}(t),t)$, is given by
\begin{align}
    \bm{v}(\bm{q}_0(t),\bm{y}(t),t) :=
    \begin{bmatrix}
        \frac{2}{2+\gamma\psi_0}\bm{f}^{(1)}(\bm{q}_0(t),\bm{y}(t),t) \\
        \frac{2}{2+\gamma\psi_0}\bm{f}^{(2)}(\bm{q}_0(t),\bm{y}(t),t) \\
        0 \\
        \vdots \\
        0
    \end{bmatrix}.
\end{align}

\subsubsection{Fourth order spatial discretization}
As in Fazio's second order approximation~\reva{\cite{FazioAndJanelli2014}}, we start by using the chain rule for the derivative of $q(\xi(x),t)$ with respect to $x$ at a node with index $n$ to get
\begin{align} \label{Chain_rule:eq}
    \frac{\partial q(\xi(x),t)}{\partial x}\biggr\rvert_{x_n} = \frac{\partial q(\xi(x),t)}{\partial \xi}\biggr\rvert_{x_n} \cdot \frac{d \xi(x)}{d x}\biggr\rvert_{x_n}.
\end{align}
Fazio approximates both derivatives on the right hand side of~\eqref{Chain_rule:eq} by centered differences~\reva{\cite{FazioAndJanelli2014}}.
In contrast, we observe that $\left(\nicefrac{d\xi(x)}{dx}\right)_{x_n}$ can be obtained exactly from~\eqref{LogRule} which leads to simpler expressions when using higher order derivative approximations with wider stencils.

Since the logarithmic map~\eqref{LogRule} between grids is a bijection, we can obtain its inverse mapping
\begin{align}
    x(\xi) = -c \ln(1-\xi) &\iff \xi(x) = 1 - e^{-\nicefrac{x}{c}} 
\end{align}
and calculate its derivative at integer points
\begin{align}
    \frac{\partial \xi(x)}{\partial x}\biggr\rvert_{x_n} &= \frac{1}{c} e^{-\nicefrac{x_{n}}{c}} = \frac{1}{c}\left( 1 - \xi_{n} \right) = \frac{1}{c}\left( 1 - \frac{n}{N} \right).
\end{align}
By plugging this into~\eqref{Chain_rule:eq} we obtain
\begin{align} \label{Chain_rule:rel}
    \frac{\partial q(\xi(x),t)}{\partial x}\biggr\rvert_{x_n} = \frac{\partial q(\xi(x),t)}{\partial \xi}\biggr\rvert_{x_n} \cdot \frac{1}{c}\left( 1 - \frac{n}{N} \right),
\end{align}
connecting deriatives on the infinite grid with derivatives on the equidistant, finite reference grid.
Equation~\eqref{Chain_rule:rel} allows to define finite difference approximations for the first derivative in a bounded uniform mesh $(\xi_n)_{n=0,\ldots,N-1}$ and then transform them to the mesh in the semi-infinite domain $(x_n)_{n=0,\ldots,N-1}$.

To approximate the derivative $\nicefrac{\partial q(\xi(x),t)}{\partial\xi}\rvert_{x_n}$ of $q$ on the equidistant grid at the interior points, we use the fourth order compact finite difference by Lele~\cite[Eq. (2.1.6)]{lele1992} with $\alpha=\frac{1}{4}$, resulting in
\begin{align}
	\label{eq:pade}
    \frac{\partial q(\xi(x),t)}{\partial \xi}\biggr\rvert_{x_{n-1}} + \frac{4 \partial q(\xi(x),t)}{\partial \xi}\biggr\rvert_{x_{n}} + \frac{\partial q(\xi(x),t)}{\partial \xi}\biggr\rvert_{x_{n+1}} \approx \frac{3}{d} \left( q_{n+1}(t) - q_{n-1}(t) \right)
\end{align}
with $d := \xi_{n+1} - \xi_{n} = \frac{1}{N}$ for the interior points with indices $n \in \{1,2,\dots,N-1\}$.
Note that non-compact higher order finite differences would require additional ghost points at the left boundary which could be difficult to realize for the nonlinear boundary condition in~\eqref{eq:MRSys}.

We use~\eqref{Chain_rule:rel} to replace the derivatives in uniform space by derivatives in pseudo-space so that~\eqref{eq:pade} becomes
\begin{align}
    \label{eq:pade:pseudo-space}
    \left( \frac{Nc}{N-n+1} \right) \frac{\partial q(\xi(x),t)}{\partial x}\biggr\rvert_{n-1} +& \left( \frac{4Nc}{N-n} \right)  \frac{\partial q(\xi(x),t)}{\partial x}\biggr\rvert_{n} + \left(  \frac{Nc}{N-n-1} \right) \frac{\partial q(\xi(x),t)}{\partial x}\biggr\rvert_{n+1} \notag  \\
    &\approx \frac{3}{d} \left( q_{n+1}(t) - q_{n-1}(t) \right)
\end{align}
Like in Section~\ref{Second:order} we obtain the second derivative by applying the first order derivative twice. 
As a consequence, two different approximations for the point at the left boundary are required: one including the boundary condition obtained by the same procedure as described by Malele et al.~\cite{malele2022} and another one without the influence of the boundary condition obtained by using the fourth order upwinding scheme given by Lele~\cite{lele1992}.

The coefficients of the upwinding approximation at the boundary that includes the left boundary condition are calculated in \ref{appendix:a} and result in the scheme
\begin{align} 
	\label{eq:upwind-boundary-1}
    c\frac{\partial q(\xi(x),t)}{\partial x}\biggr\rvert_{x_0} + \frac{Nc}{N-1}\frac{\partial q(\xi(x),t)}{\partial x}\biggr\rvert_{x_1} \approx& \left( \frac{12\alpha dc-17\gamma}{18\gamma d} \right)q_0(t) + \frac{1}{2d}q_1(t) + \notag \\
    &+ \frac{1}{2d}q_2(t) - \frac{1}{18d}q_3(t) + \\
    &- \frac{2c}{3\gamma} \left(f(\bm{q}_0(t),\bm{y}(t),t) - \frac{\partial q(\xi(x),t)}{\partial t}\biggr\rvert_{x_0} \right). \notag
\end{align}
On the other hand, since the second upwinding scheme at $n = 0$ does not include the boundary condition, the coefficients of the scheme on the uniform grid are taken directly from~\cite[Equation (4.1.4)]{lele1992}, so that
\begin{align}
    \frac{\partial q(\xi(x),t)}{\partial \xi}\biggr\rvert_{x_{0}} + 3\frac{\partial q(\xi(x),t)}{\partial \xi}\biggr\rvert_{x_{1}} \approx \frac{1}{d} \left( -\frac{17}{6}q_{0}(t) + \frac{3}{2}q_{1}(t) + \frac{3}{2}q_2(t) - \frac{1}{6}q_3(t) \right).
\end{align}
We again replace the derivatives in uniform space by derivatives in pseudo-space
\begin{align}\label{eq:upwind-boundary-2}
    c \frac{\partial q(\xi(x),t)}{\partial x}\biggr\rvert_{x_0} + \left( \frac{3Nc}{N-1} \right)\frac{\partial q(\xi(x),t)}{\partial x}\biggr\rvert_{x_1} \approx \frac{1}{d} \left( -\frac{17}{6}q_{0}(t) + \frac{3}{2}q_{1}(t) + \frac{3}{2}q_2(t) - \frac{1}{6}q_3(t) \right).
\end{align}
Using equations~\eqref{eq:pade:pseudo-space},~\eqref{eq:upwind-boundary-1} and~\eqref{eq:upwind-boundary-2} we obtain the system
\begin{subequations}
\begin{align}
    M_1 \; D_x\bm{q}(t) &\approx B_1 \bm{q}(t) + K_1(\bm{q}_0(t),\bm{y}(t),t) + V_1(\dot{\bm{q}}_0(t)), \label{CFD_sys:1} \\
    M_2 \; D_x^2 \bm{q}(t) &\approx B_2 \; D_x \bm{q}(t), \label{CFD_sys:2}
\end{align}
\end{subequations}
where $M_1$, $M_2$, $B_1$, $B_2$, $K_1$ and $V_1$ are matrices and vectors defined in~\ref{appendix:b} and $D_x\bm{q}(t)$ is the spatial discretization of $\frac{\partial}{\partial x}\bm{q}(\xi(x),t)$.
Multiplying~\eqref{CFD_sys:1} by $M_1^{-1}$ yields
\begin{subequations}
\begin{align}
    D_x\bm{q}(t) &\approx M_1^{-1} \; B_1 \bm{q}(t) + M_1^{-1} \; K_1(\bm{q}_0(t),\bm{y}(t),t) + M_1^{-1} \; V_1(\dot{\bm{q}}_0(t)), \label{CFD_sys*:1} \\
    M_2 \; D_x^2 \bm{q}(t) &\approx B_2 \; D_x \bm{q}(t). \label{CFD_sys*:2}
\end{align}
\end{subequations}
Substituting the right hand side of~\eqref{CFD_sys*:1} into~\eqref{CFD_sys*:2} yields
\begin{align}
    M_2 D_x^2 \bm{q}(t) \approx B_2 \left( M_1^{-1} \; B_1 \bm{q}(t) + M_1^{-1} \; K_1(\bm{q}_0(t),\bm{y}(t),t) + M_1^{-1} \; V_1(\dot{\bm{q}}_0(t)) \right).
\end{align}
After some algebra we obtain
\begin{align*}
    M_2 \; D_x^2 \bm{q}(t) - B_2 M_1^{-1} V_1(\dot{\bm{q}}_0(t)) &\approx B_2 M_1^{-1}B_1 \bm{q}(t) + B_2 M_1^{-1} \; K_1(\bm{q}_0(t),\bm{y}(t),t).
\end{align*}
We now use the semidiscretized version $\dot{\bm{q}}(t) = D_x^2\bm{q}(t)$ of the differential equation~\eqref{MRSys1} to transform the matrix $V_1(\dot{\bm{q}}_0(t))$ so that
\begin{align}
    V_1(\dot{\bm{q}}_0(t)) = \frac{2c}{3\gamma}
    \begin{bmatrix}
        \dot{q}^{(1)}_0(t) \\
        \dot{q}^{(2)}_0(t) \\
        0 \\
        0 \\
        \vdots \\
        0
    \end{bmatrix}
    = \frac{2c}{3\gamma}
    \underbrace{
    \begin{bmatrix}
		1 & 0 & 0 & \dots \\
		0 & 1 & 0 & \ddots \\
		0 & 0 & 0 & \ddots \\
		\vdots & \ddots & \ddots & \ddots
	\end{bmatrix}}_{=:\mathbb{P}}
	\dot{\bm{q}}(t)
    = \frac{2c}{3\gamma} \; \mathbb{P} \; \dot{\bm{q}}(t).
\end{align}
We then obtain the final approximation
\begin{align}
    \dot{\bm{q}}(t) &= A_s \bm{q}(t) + \bm{v}(\bm{q}_0(t),\bm{y}(t),t)
\end{align}
to be used in~\eqref{semi-discrete:sys} where
\begin{align}
    A_s &:= \Psi^{-1} B_2 M_1^{-1} B_1, \\
    \bm{v}(\bm{q}_0(t),\bm{y}(t),t) &:= \Psi^{-1}B_2 M_1^{-1} K_1(\bm{q}_0(t),\bm{y}(t),t), \\
    \Psi &:= M_2 - \frac{2c}{3\gamma} B_2 M_1^{-1} \mathbb{P}.
\end{align}

\subsection{Time stepping methods for the semi-discrete system}
\label{subsection:time-stepping}
The semi-discrete system~\eqref{semi-discrete:sys} can be solved numerically by using a time stepping method.
In this section we briefly state the different time stepping schemes we investigate: a second order implicit method (trapezoidal rule), a fourth order Diagonally Implicit Runge-Kutta (DIRK) method and Implicit-Explicit (IMEX) Runge-Kutta methods of order two and four.

\subsubsection{Trapezoidal Rule}
The fully implicit Trapezoidal Rule applied to~\eqref{semi-discrete:sys} yields
\begin{align}
    \underbrace{\left( \mathbb{I} - \frac{\Delta t}{2} A \right)}_{=:M_{\textrm{left}}} \bm{\eta}^{k+1} = \left( \mathbb{I} + \frac{\Delta t}{2} A \right) \bm{\eta}^{k} + \frac{\Delta t}{2} \left( \bm{\omega}^{k} + \bm{\omega}^{k+1} \right),
\end{align}
where $\bm{\eta}^{k} := \bm{\eta}(t^k)$ and $\bm{\omega}^{k} := \bm{\omega}(\bm{q}_0(t^k), \bm{y}(t^k), t^k)$.
Due to the $\bm{\omega}^{k+1}$ term on the right hand side this is an implicit system and requires a nonlinear solver.
We use the Newton-Krylov method LGMRGES with $M_{\textrm{left}}$ as preconditioner, implemented in the {\tt newton\_krylov} function of the {\tt scipy.optimize} Python library~\cite{2020SciPy-NMeth,kelley2003solving, knoll2004jacobian, baker2005technique}.
We take the solution of the explicit system
\begin{align} \label{ExpDisSys}
    \bm{\tilde{\eta}}^{k+1} = \left( \mathbb{I} + \Delta t \; A \right) \bm{\eta}^k + \Delta t \; \bm{\omega}^{k}
\end{align}
as starting value for the Newton method.

\subsection{Diagonally Implicit Runge Kutta (DIRK) method}
As implicit time integrator we use the fourth order DIRK method ESDIRK4(3)6L[2]SA~\cite{KennedyAndCarpenter2016}.
Kennedy and Carpenter recommend this method as the ``\emph{default method for solving stiff problems at moderate error tolerances}''.
In ESDIRK methods, the first stage is explicit but all other stages are implicit and require a nonlinear solver.
As for the trapezoidal rule, we use the Python function {\tt newton\_krylov}.
The previous stage is used as starting value.

\subsection{Implicit-Explicit Runge-Kutta (IMEX) methods}
In the semi-discrete system~\eqref{semi-discrete:sys}, only the $\bm{\omega}$ term arising from the boundary condition is nonlinear.
The $A$ term is linear but stiff because of the discrete Laplacian. 
To avoid the overhead of having to use a fully nonlinear solver for every stage, we propose to use an implicit-explicit RKM instead.
This allows to treat the boundary term explicitly, thus avoiding a nonlinear solver per stage, while treating the stiff discrete Laplacian implicitly, which avoids a very harsh time step restriction.
We use the second-order IMEX Midpoint rule by Ascher et al.~\cite{AscherEtAl1997} and the fourth order IMEX method combining the explicit ARK4(3)6L[2]SA–ERK and implicit ARK4(3)6L[2]SA–ESDIRK by Kennedy and Carpenter~\cite{KennedyAndCarpenter2016}.

\section{Numerical results}\label{section:results}
We compare accuracy and computational cost of (i) the finite difference (FD) schemes described in Section~\ref{section:mre},  (ii) the third-order direct numerical integrator proposed by Daitche~\cite{Daitche2013} and (iii) our reimplementation of Prasath's polynomial expansion method~\cite{Fokas1997} for the reformulated problem~\eqref{eq:MRSys} by Prasath et al.~\cite{PrasathEtAl2019}.
Note that we rely on our own implementation of Prasath's methods since no publicly available code appears to exists~\cite{UrizarnaEtAl2023}.
Daitche~\reva{\cite{Daitche2013}} provided us with his code in private communication.
We compared both our and his implementations, ensured that the coefficients computed in our version match those reported in the paper and finally reproduced figures 3 and 4 in his paper to confirm that both implementations deliver the same result.
One difference, however, is that we only compute coefficients in double instead of quadruple precision, which in some cases leads to a reduction in convergence order for the 3rd order variant.
This effect was already observed by Daitche~\cite{Daitche2013} and Moreno-Casas and Bombardelli~\cite{MorenoCasasEtAl2016}.

\subsection{Benchmark problems}
We consider three flow fields where analytical solutions to the MRGE are available: a steady vortex~\cite{CandelierEtAl2004}, a quiescent flow and an unsteady but spatially homogeneous oscillatory background~\cite{PrasathEtAl2019}.
Additionally, we consider two unsteady and inhomogeneous flow fields  where no analytic solution is known, the Bickley jet proposed by Rypina et al.~\cite{rypina2007lagrangian} with the parameters from Hadjighasem et al.~\cite{hadjighasem2016spectral} and an experimentally measured Faraday flow~\cite{colombi2022coexistence, colombi2021three}.
In the last two cases, we measure the error against a high resolution reference computed with Prasath et al.'s~\reva{\cite{PrasathEtAl2019}} algorithm.
For each flow field, we investigate the accuracy of the methods when simulating trajectories of particles that are lighter ($R < 1$) or denser ($R > 1$) than the fluid or neutrally buoyant ($R=1$) and cases with zero and non-zero initial relative velocity.
Specifically, we use $R=7/9$ where $\beta = 2/3$ and $R=4/3$ where $\beta = 3/2$.

\subsubsection{Quiescent, steady inhomogeneous and unsteady homogeneous flow field}
\label{section:analytical_solutions}
Below, we describe in detail the three benchmark problems for which an analytical solution is available.

\begin{figure}[t]
    \centering
    \includegraphics[scale=1]{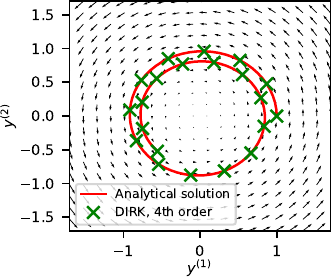}
    \includegraphics[scale=1]{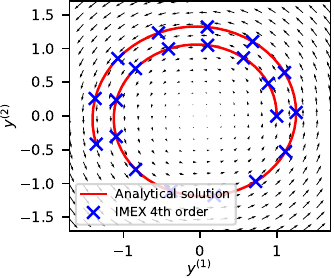}
    \caption{Trajectories of two particles moving in the vortex flow field with initial position $\bm{y}(t_0) = (1,0)^T$ and zero initial relative velocity with $S=0.3$ and $R=\nicefrac{7}{9}$ (left) or $R=\nicefrac{4}{3}$ (right) for $t\in[0,10]$.}
    \label{fig:analytic:vortex}
\end{figure}

\paragraph{Vortex flow field} Figure~\ref{fig:analytic:vortex} shows the vortex flow field and two example trajectories for different values of $R$.
The flow field
\begin{align}
    \bm{u} = \left\lVert \bm{y}(t) \right\lVert \omega \; \bm{e}_{\theta} = \omega \begin{bmatrix}
        -y^{(2)} \\
        y^{(1)}
    \end{bmatrix},
    \label{eq:vortex}
\end{align}
is a stationary vortex whose tangential velocity increases proportionally to the distance from the origin.
We set the angular velocity of the vortex to $\omega = 1$.
An analytical solution for this problem is given by Candelier et al.~\cite[Eq. (12)]{CandelierEtAl2004}. The markers in the figure indicate the numerical solution computed with a fully implicit fourth-order DIRK method (left) and a fourth-order implicit-explicit Runge-Kutta method (right). 
The analytical solutions are shown as red lines.
For both the lighter and heavier particle, the solution agrees very well with the numerical approximation.
Solutions computed with the other numerical methods are omitted to keep the plot readable but can be generated using the code provided with the paper.
For the vortex flow field, an analytical solution is only available if the initial relative velocity of the particle is zero and thus we only consider this case.

\begin{figure}[t]
    \centering
    \includegraphics[scale=1]{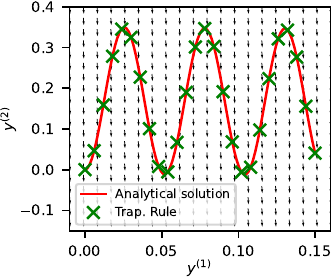}
    \includegraphics[scale=1]{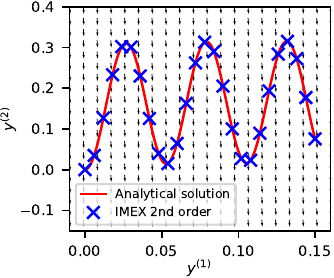}
    \caption{Trajectories of two particles moving in the oscillatory flow field with initial position $\bm{y}(t_0) = {(0,0)}^T$ and zero initial relative velocity with $S=0.3$ and $R=\nicefrac{7}{9}$ (left) or $R=\nicefrac{4}{3}$ (right) for $t\in[0,3]$. The arrows show the velocity field at $t=3$.}
    \label{fig:analytic:oscillatory}
\end{figure}

\paragraph{Oscillatory flow field} Figure~\ref{fig:analytic:oscillatory} shows the analytical and numerical solutions for particles that are lighter (left) and denser (right) than the fluid in a spatially homogeneous oscillatory flow field 
\begin{align}
    \bm{u}(t) =
    \begin{bmatrix}
        u^{(1)} \\
        \sin(\lambda t)
    \end{bmatrix},
\end{align}
where $u^{(1)} = 0.05 $ is the horizontal component of the flow field and $\lambda = 6$ is the oscillation frequency.
We show the numerical approximations provided by the trapezoidal rule (left) and the second-order IMEX method (right) as markers. 
Again, the solutions for the other numerical methods can be generated using the provided code.
The analytical solution reads
\begin{subequations}
\begin{align}
    \label{eq:analytical:quiescent}
    y^{(1)}(t) =& \; y^{(1)}_0 + u^{(1)} t + 2\frac{v_0^{(1)}}{\pi}  \int_{0}^{\infty} \frac{\gamma(1- e^{-k^2t})}{(\alpha -k^2)^2+(k\gamma)^2} dk, \\
    y^{(2)}(t) =& \; y_0^{(2)} + \frac{1-\cos(\lambda t)}{\lambda} + \frac{2}{\pi} \int_0^\infty \; v_0^{(2)} \frac{ \gamma (1 - e^{-k^2 t}) }{k^2\gamma^2 + (k^2 - \alpha)^2} \; dk + \notag \\
    &+ \frac{2}{\pi} \frac{(1-R)\lambda}{R} \int_0^\infty \frac{k^2\gamma e^{-k^2 t}}{(k^2\gamma^2+(k^2-\alpha)^2)(k^4 + \lambda^2)}dk + \notag \\
    &+ \frac{2}{\pi} \frac{(1-R)\lambda}{R} \int_0^\infty \frac{k^4 \gamma \sin(\lambda t)}{(k^2\gamma^2+(k^2-\alpha)^2)(k^4 + \lambda^2)\lambda} dk + \notag \\
    \label{eq:analytical:oscillatory}
    &- \frac{2}{\pi} \frac{(1-R)\lambda}{R} \int_0^\infty \frac{k^2\gamma \cos(\lambda t) }{(k^2\gamma^2+(k^2-\alpha)^2)(k^4 + \lambda^2)}dk.
\end{align}
\end{subequations}
It is shown as a solid red line in Figure~\ref{fig:analytic:oscillatory}.
The expressions for the horizontal~\eqref{eq:analytical:quiescent} and vertical~\eqref{eq:analytical:oscillatory} components are obtained by using $f^{(1)}\equiv 0$ and $f^{(2)}(s) = \left(\frac{1}{R}-1\right)\left(\lambda \cos(\lambda s)\right)$ in equation (3.16) in Prasath et al.~\cite{PrasathEtAl2019}.
As for the vortex, the numerical approximations shown as crosses match the analytical solutions.

\begin{figure}[t]
    \centering
    \includegraphics[scale=1]{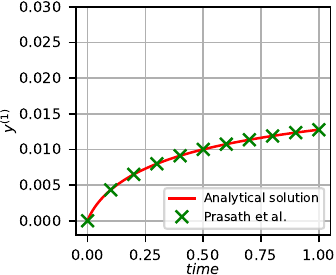}
    \includegraphics[scale=1]{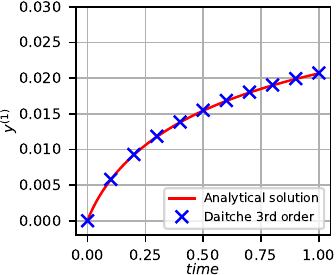}
    \caption{Horizontal component of the position versus time for two particles decelerating in a quiescent flow with initial position $\bm{y}(t_0) = {(0,0)}^T$ and initial relative velocity $\bm{q}(0, t_0) = {(0.1, 0)}^T$ with $S=0.3$ and $R=\nicefrac{7}{9}$ (left) or $R=\nicefrac{4}{3}$ (right)  for $t\in[0,1]$.}
    \label{fig:analytic:quiescent}
\end{figure}

\paragraph{Quiescent flow field} 
In this example, the flow field is zero everywhere and the particle, starting with a non-zero initial velocity, decelerates until it comes to a halt.
The solution for a single component corresponds to equation~\eqref{eq:analytical:quiescent} with $u^{(1)}=0$.
Figure~\ref{fig:analytic:quiescent} again shows that numerical and analytical solution agree.
Since no motion is happening if the particle is initiall\reva{y} at rest, we only consider the case of non-zero initial relative velocity.

\subsubsection{Fully unsteady and inhomogeneous flow fields}
Next we describe the two benchmark problems with flow fields that vary in both space and time where no analytical solution is available.

\begin{figure}[t]
    \centering
    \includegraphics[scale=1]{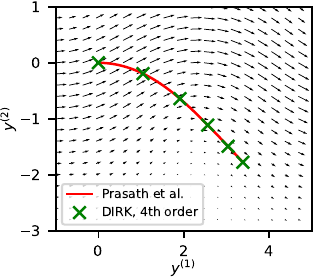}
    \includegraphics[scale=1]{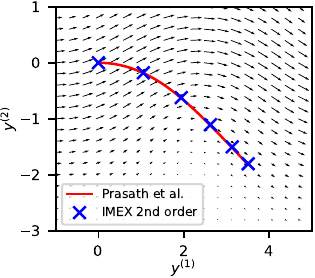}
    \caption{Trajectory of two particles moving in the Bickley jet with initial position $\bm{y}(t_0) = {(0, 0)}^T$ and  zero initial relative velocity with $S=0.3$ and $R=\nicefrac{7}{9}$ (left) and $R=\nicefrac{4}{3}$ (right) for $t\in[0,1]$. The arrows show the velocity field at $t=1$.}
    \label{fig:analytic:bickley}
\end{figure}

\begin{figure}[t]
    \centering
    \includegraphics[scale=1]{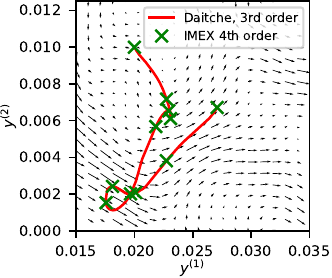}
    \includegraphics[scale=1]{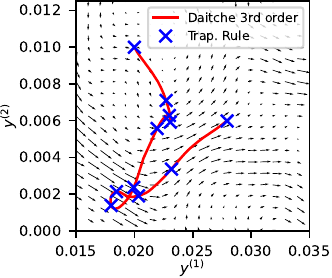}
    \caption{Trajectory of two particles moving in the Faraday flow with initial position $\bm{y}(t_0) = {(0.02, 0.01)}^T$ and zero initial relative velocity, i.e. $\bm{q}(0, t_0) = {(0, 0)}^T$, with $S=0.3$ and either $R=\nicefrac{7}{9}$ (left) or $R=\nicefrac{4}{3}$ (right) for $t\in[0,5]$. The arrows show the velocity field at $t=5$.}
    \label{fig:analytic:faraday}
\end{figure}

\paragraph{Bickley jet}
The Bickley jet is a flow field with stream function
\begin{equation}
    \Psi(x,y,t) = -U_0 L \text{tanh}(y/L) + \sum_{i=1}^3 A_i U_0 L \text{sech}^2(y/L) \cos(k_i x - \sigma_i t)
\end{equation}
and was proposed by Rypina et al.~\cite{rypina2007lagrangian} as a model for atmospheric jets.
We use the values for $U_0$, $L$, $k_i$ and $\sigma_i$ given by Hadjighasem et al.~\cite{hadjighasem2016spectral}. 
This flow field ``serves as an idealized model of the stratospheric flow"~\cite{padberg2017network}.
Figure~\ref{fig:analytic:bickley} shows the reference solution computed with Prasath et al.'s algorithm~\reva{\cite{PrasathEtAl2019}} and solutions computed with the second- and fourth-order IMEX/finite difference approach.
\paragraph{Faraday flow}
The second unsteady, inhomogeneous  field is an experimentally measured Faraday flow, arising in a thin layer of water in a circular container attached to a shaker vibrating at a frequeny of \SI{50}{\hertz}. 
Details of the experiment are given by Colombi et al.~\cite{colombi2022coexistence, colombi2021three}.
Measurements are taken at $115 \times 86$ points in space with \SI{2}{\milli\meter} resolution with a snapshot in time taken every \SI{40}{\milli\second} and a total of \SI{1056} snapshots.
We use rectangular bivariate splines provided by \texttt{scipy} to interpolate between data points in space to obtain the velocity field at a given particle position and linear interpolation to interpolate between snapshots in time.
Figure~\ref{fig:analytic:faraday} shows example trajectories of particles in the Faraday flow computed with two different algorithms for particles that are lighter (left) and heavier (right) than the fluid.
\subsection{Convergence order and accuracy}
This section analyzes the convergence orders of all three methods, Prasath al., Daitche and FD, for the five flow fields described above. 
We investigate the impact of changing the values of $R$ and $S$ and of whether the particle has zero or non-zero initial relative velocity.
For the flow fields where an analytic solution is available we report the maximum relative $l_2$-error over all time steps.
If no analytical solution is known and a numerically computed reference is used, we report the relative $l_2$-error at the final time step.
We choose values of $R=\nicefrac{1}{3}$, corresponding to a particle of gas with density $\rho_p \to 0$ and $R=\nicefrac{7}{3}$, corresponding to the density of a particle of a heavy rock such as Basalt~\cite{sharma1997environmental}.
In addition, we simulate three intermediate cases of a light particle with $R=\nicefrac{7}{9}$, a neutrally buoyant particle with $R=1$ and a heavy particle with $R=\nicefrac{4}{3}$. 
Values for the Stokes number are $S\in\{0.01, 0.1, 0.5, 1, 2, 4\}$.

\subsubsection{Quiescent, steady inhomogeneous and unsteady homogeneous flow field}
We first analyze convergence for the three flow fields for which we have analytical solution.
Throughout, $N$ denotes the number of nodes in space and steps in time for the FD and Prasath et al.'s~\reva{\cite{PrasathEtAl2019}} polynomial expansion method, which we always keep equal.
For Daitche's method~\reva{\cite{Daitche2013}}, where no spatial discretization is needed, $N$ simply indicates the number of time steps.

\begin{figure}[t]
    \centering
    \includegraphics[scale=1]{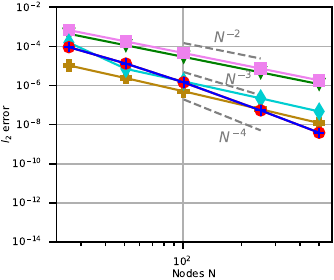}
    \includegraphics[scale=1]{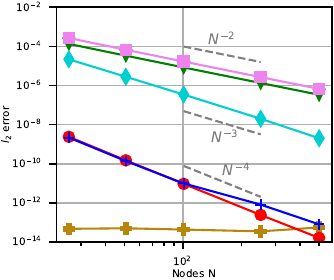}
    \includegraphics[scale=1]{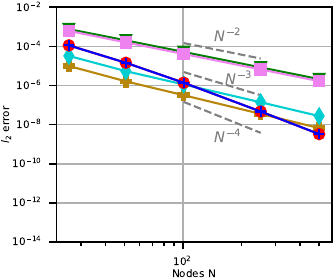}
    \includegraphics[scale=1]{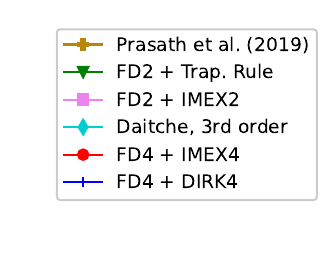}
    \caption{Error versus $N$ for a particle moving in the vortex with initial position $\bm{y}(t_0) = {(1,0)}^T$ and zero initial relative velocity with $S=0.1$ and $R=\nicefrac{7}{9}$ (top left), $R=1$ (top right) or $R=\nicefrac{4}{3}$ (bottom left) for $t\in[0,1]$.}
    \label{fig:convergence:vortex}
\end{figure}

\paragraph{Vortex flow field}
Figure~\ref{fig:convergence:vortex} shows convergence of all methods against the analytical solution for three different values of $R$. 
Lines with slopes two, three and four are shown as a guide to the eye.
Second order FD converge with their theoretically expected rate in all three cases, both with a diagonally implicit or IMEX Runge-Kutta method.
Fourth-order FD with DIRK or IMEX converges with orders between 3 and 4, slightly lower than what would be theoretically expected unless, for $R \neq 1$ and order 4 if $R=1$.
Daitche's~\reva{\cite{Daitche2013}} third order method achieves its expected convergence order only for the case of a neutrally buoyant particle ($R=1$, middle) and converges with about order two for the lighter or denser particle.
Most likely this is because coefficients were computed in double and not quadruple precision~\cite{MorenoCasasEtAl2016,Daitche2013}.
For very small values of $S$ and $R$, Daitche's~\reva{\cite{Daitche2013}} method becomes unstable.
Prasath's method converges with order between two and three for $R \neq 1$ and achieves what looks like spectral accuracy for $R=1$, where it reaches round-off error with only very few nodes.

\paragraph{Oscillatory background flow field}

\begin{figure}[t]
    \centering
    \includegraphics[scale=1]{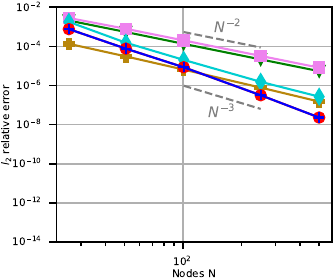}
    \includegraphics[scale=1]{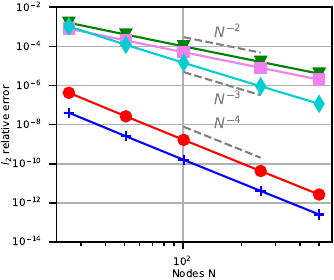}
    \includegraphics[scale=1]{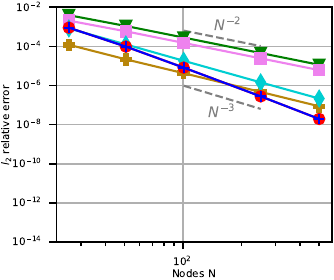}
    \includegraphics[scale=1]{Figure_legend.pdf}
    \caption{Convergence for a particle moving in an ocillatory background starting from $\bm{y}(t_0) = {(0,0)}^T$ and zero initial relative velocity, $S=0.1$ and either $R=\nicefrac{7}{9}$ (top left), $R=1$ (top right) or $R=\nicefrac{4}{3}$ (bottom left) for $t\in[0,1]$}
    \label{fig:convergence:oscillatory1}
\end{figure}

\begin{figure}[t]
    \centering
    \includegraphics[scale=1]{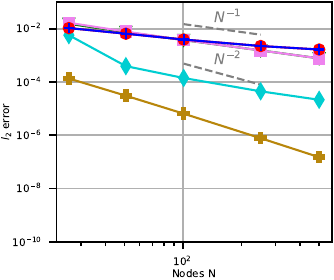}
    \includegraphics[scale=1]{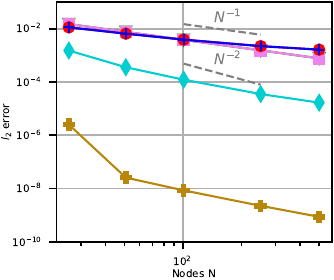}
    \includegraphics[scale=1]{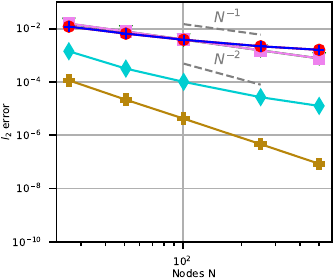}
    \includegraphics[scale=1]{Figure_legend.pdf}
    \caption{Convergence for a particle moving in an oscillatory background starting from $\bm{y}(t_0)={(0,0)}^T$ and initial relative velocity $\bm{q}(0, t_0)={(0,0.1)}^T$, with $S=0.1$ and $R=\nicefrac{7}{9}$ (top left), $R=1$ (top right) or $R=\nicefrac{4}{3}$ (bottom left) for $t\in[0,1]$.}
    \label{fig:convergence:oscillatory2}
\end{figure}

Figures~\ref{fig:convergence:oscillatory1} shows convergence of all methods for particles of different densities with zero relative velocity in the oscillatory flow field.
Figure~\ref{fig:convergence:oscillatory2} shows convergence for the same particles but with a small non-zero initial relative velocity.
For zero initial relative velocity, results are very similar to those for the vortex. 
Second-order FD converges with its expected order with both implicit trapezoidal rule and IMEX2.
Fourth-order FD reaches its full order only for $R=1$ and converges with about order three otherwise.
Daitche~\reva{\cite{Daitche2013}} achieves order three for $R=1$ and with order somehwat higher than two if $R \neq 1$.
As before, Prasath converges with order between two and three unless the particle is neutrally buoyant, in which case it again achieves spectral convergence.

The picture changes dramatically if the particle starts with a non-zero initial relative velocity.
Figure~\ref{fig:convergence:oscillatory2} shows convergence for the same parameters as above except that now the velocity $\bm{v}(t_0)$ of the particle at the beginning is slightly different from the flow field at its starting point, leading to a small but non-zero initial relative velocity of $\bm{q}(0,t_0) = (0, 0.1)$.
This leads to a discontinuity in the initial value for~\eqref{eq:MRSys} since we have $q^{(2)}(0,t_0) > 0$ but $q^{(2)}(x,t_0) = 0$ for $x > 0$.
For FDM, this leads to a massive deterioration in convergence order.
Second-order FD methods are reduced to firsrt order whereas the fourth-order FD variants converge even slower, with orders of about 0.7.
Prasath's method looses the spectral accuracy we observed for zero initial relative velocity but still remains highly accurate and converges with approximately order two.
While finite difference rely on Taylor series and thus are well-known to have big problems resolving discontinuities, Prasath's method uses an integral formulation which seems to be less affected.
However, even Daitche's method~\reva{\cite{Daitche2013}}, which does not rely on~\eqref{eq:MRSys}, is also reduced to slightly better than first order convergence.
\reva{One reason is likely the very stiff initial relaxation of the particle velocity to the flow field.
Order reduction for stiff ODEs is a well-documented phenomenon~\cite{FrankEtAl1985}.
However, the choice of the parameter $c$, see \S\ref{subsec:fd}, and other factors probably also play a role.}

\begin{figure}[th]
    \centering
    \includegraphics[scale=1]{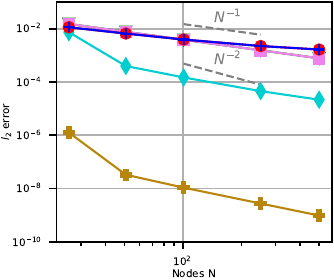}
    \includegraphics[scale=1]{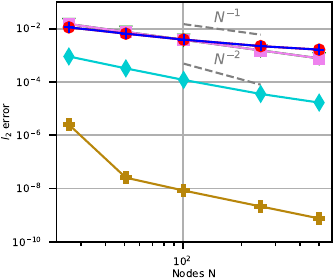}
    \includegraphics[scale=1]{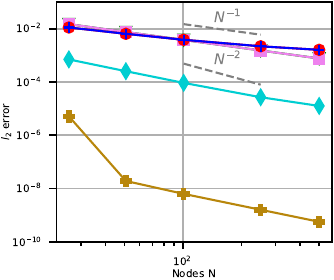}
    \includegraphics[scale=1]{Figure_legend.pdf}
    \caption{Convergence plot for a particle moving in a quiescent flow, initialised at $\bm{y}(t_0)={(0,0)}^T$ and initial relative velocity $\bm{q}(0,t_0)={(0.1,0)}^T$, with $S=0.1$ and $R=\nicefrac{7}{9}$ (top left), $R=1$ (top right) or $R=\nicefrac{4}{3}$ (bottom left) for $t\in[0,1]$.}
    \label{fig:convergence:quiescent}
\end{figure}

\paragraph{Relaxing particle}
Figure~\ref{fig:convergence:quiescent} shows the convergence plots for all methods against the analytical solution of a relaxing particle in a quiescent flow~\eqref{eq:analytical:quiescent} for three different values of $R$. 
Since the field is zero everywhere, the particle has to start with nonzero initial velocity to induce motion.
The result is very similar to what we showed for the oscillatory background flow field.
FD2 converges with order one with FD4 converging even slower.
Daitche~\reva{\cite{Daitche2013}} converges slightly faster than order one but is more accurate while Prasath is by far the most accurate and converges with order two or more.

\subsubsection{Fully unsteady and inhomogeneous flow fields}
Here we study convergence for the two unsteady and inhomogeneous fields where the error is measured against a numerically computed reference solution.

\paragraph{Bickley Jet}
Tables~\ref{table:convergence:bickley1} and~\ref{table:convergence:bickley2} show convergence rates of all methods for a particle in the Bickley jet with different values of $R$ and $S$ with zero and non-zero initial relative velocity.
We observe a very similar pattern as for the three simpler test cases before: for zero initial relative velocity, FD2 and Daitche~\reva{\cite{Daitche2013}} largely achieve close to their theoretical order.
For very small values of $S$ and $R$, however, Daitche's method~\reva{\cite{Daitche2013}} eventually becomes unstable.
The fourth order FD method suffers from order reduction for very light or very heavy particles, although the effect is more pronounced in the former case, as well as for very small Stokes numbers, but still provides the highest order in almost all cases.
In these cases, relaxation of the particle velocity against the flow field is very rapid and the problem becomes very stiff.
Prasathet al.'s method converges with orders between two and four, depending on the values of $R$ and $S$.
A notable exception is zero initial relative velocity and $R=1$ where it achieves machine precision.

For nonzero initial relative velocity, the finite difference method suffers from dramatic order reduction.
FD2 at least remains first order accurate throughout but FD4 is reduced to about order $0.5$, rendering it largely useless in this case.
Most likely this is because the initial value in~\eqref{eq:MRSys} is not continuous for a non-zero initial relative velocity and finite difference are well known to cope badly with discontinuities.
Daitche's method~\reva{\cite{Daitche2013}} also suffers from reduced order but not as dramatic as the FDM, achieving orders between one and two.
By contrast, convergence of Prasath's method remains largely unaffected, probably because it relies on an integral formulation of the PDE which can better handle discontinuities.
It delivers the highest order in all cases.

\begin{table}[t]
    \begin{center}
    \begin{tabular}{|| c | c | c | c | c | c ||} 
        \hline
        \multicolumn{6}{|c|}{Bickley jet with zero initial relative velocity and $S=0.1$} \\        
        \hline
        Method & $R=1/3$ & $R=7/9$ & $R=1$ & $R=4/3$ & $R=7/3$ \\ 
        \hline\hline
        Prasath et al.  & $2.3$ & $2.5$ &  \textcolor{red}{machine-$\varepsilon$} & $2.6$ & $3.2$ \\
        \hline
        FD2 + Trap. & $2.1$ & $2.1$ & $2.0$ & $2.0$ & $2.00$ \\
        \hline
        FD2 + IMEX2 & $2.1$ & $2.1$ & $2.0$ & $2.0$ & $2.0$ \\
        \hline
        Daitche3  & \textcolor{red}{unstable} & $2.9$ & $3.0$ & $2.8$ & $2.7$\\
        \hline
        FD4 + IMEX4 & $\textbf{2.8}$ & $\textbf{3.2}$ & $\textbf{4.1}$ & $\textbf{3.2}$ & $\textbf{3.5}$ \\
        \hline
        FD4 + DIRK4 & $\textbf{2.8}$ & $\textbf{3.2}$ & $3.9$ & $\textbf{3.2}$ & $\textbf{3.5}$ \\
        \hline
    \end{tabular}
    \end{center}
    %%
    %\newline
    \vspace{0em}
    %\newline
    %%
    \begin{center}
    \begin{tabular}{|| c | c | c | c | c | c | c ||}
        \hline
        \multicolumn{7}{|c|}{Bickley jet with zero initial relative velocity and $R=\nicefrac{4}{3}$} \\            
        \hline
        Method & $S=0.01$ & $S=0.1$ & $S=0.5$ & $S=1$ & $S=2$ & $S=4$ \\ [0.5ex] 
        \hline\hline
        Prasath et al. & $2.4$ & $2.6$ & $2.7$ & $3.4$ & $3.0$ & $\textbf{4.2}$ \\
        \hline
        FD2 + Trap. & $2.0$ & $2.0$ & $2.0$ & $2.0$ & $2.0$ & $2.0$ \\
        \hline
        FD2 + IMEX2 & $2.0$ & $2.0$ & $2.0$ & $2.0$ & $2.0$ & $2.0$ \\
        \hline
        Daitche3  & \textcolor{red}{unstable} & $2.8$ & $2.8$ & $2.9$ & $2.9$ & $3.0$ \\
        \hline
        FD4 + IMEX4 & $\textbf{2.4}$ & $\textbf{3.2}$ & $\textbf{4.0}$ & $\textbf{4.0}$ & $\textbf{3.9}$ & $3.9$ \\ \hline
        FD4 + DIRK4 & $\textbf{2.4}$ & $\textbf{3.2}$ & $\textbf{4.0}$ & $\textbf{4.0}$ & $\textbf{3.9}$ & $3.8$ \\        
        \hline
    \end{tabular}    
    \end{center}
    \caption{Numerically computed convergence orders for a particle in the Bickley jet with for $t\in[0,1]$ and different values of $R$ and $S$. The initial conditions are $\bm{y}={(0,0)}^T$ and $\bm{q}(0,t_0)={(0,0)}^T$.}\label{table:convergence:bickley1}
\end{table}

\begin{table}[thp]
    \begin{center}
    \begin{tabular}{|| c | c | c | c | c | c ||} 
        \hline
        \multicolumn{6}{|c|}{Bickley jet with non-zero initial relative velocity and $S=0.1$} \\    
        \hline
        Method & $R=1/3$ & $R=7/9$ & $R=1$ & $R=4/3$ & $R=7/3$ \\ 
        \hline\hline
        Prasath et al. & $\textbf{2.3}$ & $\textbf{2.5}$ & $\textbf{2.6}$ & $\textbf{2.7}$ & $\textbf{3.3}$ \\
        \hline
        FD2 + Trap./IMEX2 & $0.9$ & $0.9$ & $1.0$ & $1.0$ & $1.0$ \\
        \hline
        Daitche3 & \textcolor{red}{unstable} & $1.9$ & $1.4$ & $1.3$ & $1.3$ \\
        \hline
        FD4 + IMEX4 or DIRK4& $0.5$ & $0.5$ & $0.5$ & $0.5$ & $0.6$ \\
        \hline
    \end{tabular}
    \end{center}
    %
    %\newline
    \vspace*{0em}
    %\newline
    \begin{center}
        \begin{tabular}{|| c | c | c | c | c | c | c ||}
        \hline
        \multicolumn{7}{|c|}{Bickley jet with non-zero initial relative velocity and $R=\nicefrac{4}{3}$} \\            
        \hline
        Method & $S=0.01$ & $S=0.1$ & $S=0.5$ & $S=1$ & $S=2$ & $S=4$ \\ [0.5ex] 
        \hline\hline
        Prasath et al. & $\textbf{2.3}$ & $\textbf{2.7}$ & $\textbf{3.0}$ & $\textbf{3.5}$ & $\textbf{3.4}$ & $\textbf{4.1}$ \\
        \hline
        FD2 + Trap./IMEX2 & $1.00$ & $1.0$ & $1.0$ & $1.0$ & $1.0$ & $1.0$ \\
        \hline
        Daitche3 & \textcolor{red}{unstable} & $1.3$ & $1.3$ & $1.4$ & $1.5$ & $1.5$ \\
        \hline
        FD4 + IMEX4 & $0.5$ & $0.5$ & $0.6$ & $0.6$ & $0.6$ & $0.7$ \\ \hline
        FD4 + DIRK4 & $0.5$ & $0.5$ & $0.6$ & $0.6$ & $0.6$ & $0.7$ \\         
        \hline
    \end{tabular}
    \end{center}
    \caption{Numerically computed convergence orders for a particle in the Bickley jet for $t\in[0,1]$ with different values for $S$ and $R$. The initial conditions are $\bm{y}={(0,0)}^T$ and $\bm{q}(0,t_0)={(0.5414,0)}^T$.}
    \label{table:convergence:bickley2}
\end{table}

\paragraph{Faraday flow}
Tables~\ref{table:convergence:faraday1} and~\ref{table:convergence:faraday2} show convergence rates of all methods for a particle in the Faraday flow field with different values of $R$ and $S$ and zero and non-zero initial relative velocity.
The results are similar to those for the Bicklet jet, although Daitche~\reva{\cite{Daitche2013}} and FDM show somewhat lower convergence orders, even for zero initial relative velocity.
Second-order FD is second order for zero and about first order for non-zero initial relative velocity, although FD-IMEX2 is somewhat below order one for large values of $S$.
The instability of Daitche's method~\reva{\cite{Daitche2013}} for very small values of $S$ and $R$ as well as the dramatic order reduction for fourth-order FD for non-zero initial relative velocity remain: in the latter case, FD4 converges only with order slightly better than one-half.
Prasath et al.'s~\reva{\cite{PrasathEtAl2019}} method provides the highest order in almost all cases with non-zero initital relative velocity.
In contrast to the Bickley jet, there is some variation in which method provides the highest order for zero relative initial velocity.
Depending on the specific values for $S$ and $R$, it can be Prasath~\reva{\cite{PrasathEtAl2019}}, Daitche~\reva{\cite{Daitche2013}} or FD4.

\begin{table}[th]
    \begin{center}
    \begin{tabular}{|| c | c | c | c | c | c ||} 
        \hline
        \multicolumn{6}{|c|}{Faraday flow with zero initial relative velocity and $S=0.1$} \\    
        \hline
        Method & $R=1/3$ & $R=7/9$ & $R=1$ & $R=4/3$ & $R=7/3$ \\ 
        \hline\hline
        Prasath et al. & $\textbf{2.2}$ & $2.5$ & $2.9$ & $2.6$ & $\textbf{2.8}$ \\
        \hline
        FD2 + Trap. & $2.2$ & $2.1$ & $2.1$ & $1.9$ & $2.2$ \\
        \hline
        FD2 + IMEX2 & $2.1$ & $2.1$ & $2.1$ & $2.3$ & $2.1$ \\
        \hline
        Daitche3 & \textcolor{red}{unstable} & $\textbf{3.2}$ & $2.7$ & $2.6$ & $2.5$ \\
        \hline
        FD4 + IMEX4 & $2.1$ & $2.6$ & $\textbf{4.0}$ & $\textbf{2.7}$ & $2.7$ \\
        \hline
        FD4 + DIRK4 & $2.1$ & $2.0$ & $2.1$ & $2.6$ & $2.6$\\
        \hline
    \end{tabular}
    \end{center}
    \vspace*{0em}
    \begin{center}
        \begin{tabular}{|| c | c | c | c | c | c | c ||}
        \hline
        \multicolumn{7}{|c|}{Faraday flow with zero initial relative velocity and $R=\nicefrac{4}{3}$} \\    
        \hline
        Method & $S=0.01$ & $S=0.1$ & $S=0.5$ & $S=1$ & $S=2$ & $S=4$ \\ [0.5ex] 
        \hline\hline
        Prasath et al. & $2.3$ & $2.6$ & $2.4$ & $\textbf{2.7}$ & $2.6$ & $\textbf{3.1}$ \\
        \hline
        FD2 + Trap. & $2.1$ & $1.9$ & $1.9$ & $2.1$ & $2.1$ & $2.2$ \\
        \hline
        FD2 + IMEX2 & $2.0$ & $2.3$ & $2.1$ & $2.0$ & $2.0$ & $2.0$ \\
        \hline
        Daitche3 & \textcolor{red}{unstable} & $2.6$ & $\textbf{2.5}$ & $2.4$ & $2.4$ & $2.5$ \\
        \hline
        FD4 + IMEX4 & $2.3$ & $\textbf{2.7}$ & $2.5$ & $2.5$ & $2.2$ & $2.7$ \\
        \hline
        FD4 + DIRK4 & $\textbf{2.3}$ & $2.6$ & $2.2$ & $2.3$ & $\textbf{2.7}$ & $2.8$ \\
        \hline
    \end{tabular}
    \end{center}
    \caption{Numerically computed convergence orders for a particle in the Faraday flow field for $t\in[0,1]$ and different values of $S$ and $R$. The initial conditions are $\bm{y}={(0.02,0.01)}^T$ and $\bm{q}(0,t_0)={(0,0)}^T$.}
    \label{table:convergence:faraday1}    
\end{table}

\begin{table}[th]
    \begin{center}
    \begin{tabular}{|| c | c | c | c | c | c ||} 
        \hline
        \multicolumn{6}{|c|}{Faraday flow with non-zero initial relative velocity and $S=0.1$} \\    
        \hline
        Method & $R=1/3$ & $R=7/9$ & $R=1$ & $R=4/3$ & $R=7/3$ \\ 
        \hline\hline
        Prasath et al. & $\textbf{2.2}$ & $2.5$ & $\textbf{2.6}$ & $\textbf{2.7}$ & $\textbf{2.8}$ \\
        \hline
        FD2 + Trap. & $0.9$ & $1.0$ & $1.0$ & $1.0$ & $1.0$ \\
        \hline
        FD2 + IMEX2 & $1.1$ & $1.0$ & $1.0$ & $1.0$ & $0.9$ \\
        \hline
        Daitche3 & \textcolor{red}{unstable} & $\textbf{2.7}$ & $2.1$ & $2.2$ & $2.5$\\
        \hline
        FD4 + IMEX4 & $0.6$ & $0.6$ & $0.6$ & $0.6$ & $0.6$ \\
        \hline
        FD4 + DIRK4 & $0.6$ & $0.6$ & $0.6$ & $0.6$ & $0.6$ \\
        \hline
    \end{tabular}
    \end{center}
    \vspace*{0em}
    \begin{center}
     \begin{tabular}{|| c | c | c | c | c | c | c ||}
        \hline
        \multicolumn{7}{|c|}{Faraday flow with non-zero initial relative velocity and $R=\nicefrac{4}{3}$} \\        
        \hline
        Method & $S=0.01$ & $S=0.1$ & $S=0.5$ & $S=1$ & $S=2$ & $S=4$ \\
        \hline\hline
        Prasath et al. & $\textbf{2.3}$ & $\textbf{2.7}$ & $\textbf{2.6}$ & $\textbf{3.1}$ & $\textbf{2.6}$ & $\textbf{3.4}$ \\
        \hline
        FD2 + Trap. & $0.9$ & $1.0$ & $1.0$ & $1.0$ & $1.0$ & $1.0$ \\
        \hline
        FD2 + IMEX2 & $1.0$ & $1.0$ & $0.9$ & $0.9$ & $0.8$ & $0.7$ \\
        \hline
        Daitche3 & \textcolor{red}{unstable} & $2.2$ & $2.4$ & $2.4$ & $2.4$ & $2.6$ \\
        \hline
        FD4 + IMEX4 & $0.6$ & $0.6$ & $0.6$ & $0.6$ & $0.6$ & $0.6$ \\
        \hline
        FD4 + DIRK4 & $0.6$ & $0.6$ & $0.6$ & $0.6$ & $0.6$ & $0.6$ \\
        \hline
    \end{tabular}  
    \end{center}
    \caption{Numerically computed convergence orders for a particle in the Faraday flow field with $S=0.1$ and different values of $R$ for $t\in[0,1]$. The initial conditions are $\bm{y}={(0.02,0.01)}^T$ and $\bm{q}(0,t_0)={(0.00052558,-0.00064947)}^T$.}
    \label{table:convergence:faraday2}
\end{table}

\subsection{Computational efficiency}\label{subsec:efficiency}
This last part analyzes computational efficiency of all methods by means of work-precision plots for the two unsteady, inhomogeneous flow fields.
Work-precision plots show computational effort, here measured as runtime in seconds, on the x-axis against error on the y-axis.
This allows to quickly see either which method is the fastest to produce a solution of a certain numerical accuracy or which method can produce the most accurate solution within a certain runtime.
Runtime is measured with the \texttt{time} module in Python and includes the time required to compute all time steps of the trajectory but excludes setup times.

\reva{Note that the Fortran implementation of the method by Jaganathan et al.~\cite{JaganathanEtAl2023} for the oscillatory background flow field achieved $l_2$-errors or $10^{-4}$ in $10^{-3}$ seconds wallclock time for a $5$ second trajectory.
While we only show efficiency results for Bickley jet and Faraday flow here, results with the same code published elsewhere~\cite{JulioDiss} show that our Python implementation of Daitche's method~\reva{\cite{Daitche2013}} requires around $10^{-2}$ seconds wallclock time to match this accuracy for a $1$ second trajectory.
Whether this better performance is due to the numerical method or because compiled languages like Fortran are significantly faster than interpreted languages like Python is unknown.
Refactoring the Python code into Fortran and providing a detailed comparison is beyond the scope of this paper but would be a very useful undertaking in the future.}

\begin{figure}[t]
    \centering
    \includegraphics[scale=1]{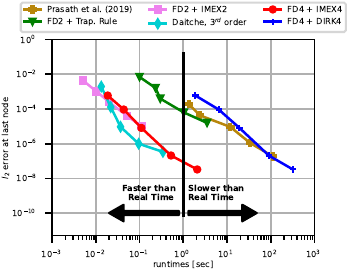}
    \includegraphics[scale=1]{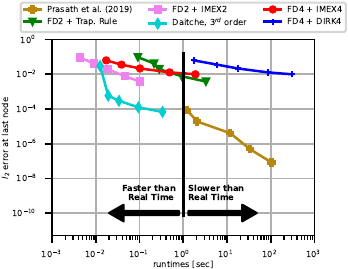}
    \caption{Error against time-to-solution for a particle moving in the Bickley Jet with $S=0.1$ and $R=\nicefrac{7}{9}$ for $t\in[0,1]$. The initial conditions are $\bm{y}(t_0) = {(0,0)}^T$ and $\bm{q}(0,t_0)={(0,0)}^T$ (left) and $\bm{q}(0,t_0)={(0.5414,0)}^T$ (right).}
    \label{fig:precision:bickley}
\end{figure}

\paragraph{Bickley jet}
Figure~\ref{fig:precision:bickley} shows error versus runtime for a particle in the Bickley jet with $S=0.1$ and $R=\nicefrac{7}{9}$ and zero (left) and non-zero (right) initial relative velocity.
The trajectory of the particle is simulated over \SI{1}{\second}, so a runtime of less than \SI{1}{\second} is considered faster than real-time even though effective use of the model for closed loop control will not be realistically possible throughout that whole range.
For zero initial relative velocity, Daitche~\reva{\cite{Daitche2013}}, FD2 + IMEX2 and FD4 + IMEX4 can all provide errors below $10^{-2}$ with runtimes between \SI{10}{\milli\second} to \SI{100}{\milli\second}.
Prasath et al.'s~\reva{\cite{PrasathEtAl2019}} method is too computationally costly to be capable of producing real-time results.
For non-zero initial relative velocity, the performance of the FD methods suffers significantly from the reduction to order one.
While FD2 + IMEX2 can still produce errors of almost $10^{-2}$ with runtimes of around \SI{100}{\milli\second}, Daitche's method~\reva{\cite{Daitche2013}} is clearly superior in this case.
Note that the provided code is written in Python and is not using multi-threading.
An optimized implementation in a compiled language would be several times faster.

\begin{figure}[t]
    \centering
    \includegraphics[scale=1]{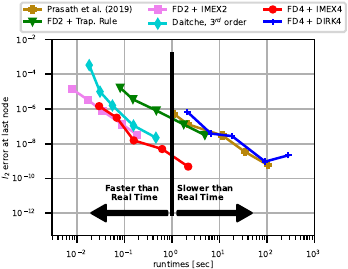}
    \includegraphics[scale=1]{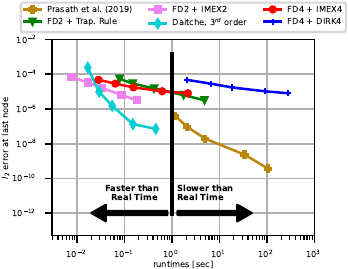}
    \caption{Error against time-to-solution for a particle moving in the Faraday flow field with $S=0.1$ and $R=\nicefrac{7}{9}$ for $t\in[0,1]$. The initial position is $\bm{y}(t_0) = {(0.02,0.01)^T}$ and the initial relative velocity is $\bm{q}(t_0)={(0,0)}^T$ (left) and $\bm{q}(0,t_0)={(0.00052558, -0.00064947)}^T$ (right).}
    \label{fig:precision:faraday}
\end{figure}

\reva{Finally, when comparing the fully implicit FD2 + Trapezoidal rule and FD4 + DIRK4 methods with their implicit-explicit counterparts FD2 + IMEX2 and FD4 + IMEX4, both show similar convergence behavior.
However, the IMEX variants require noticeably less runtime, making them more efficient.
This underscores the value of the approach in avoiding a costly fully nonlinear solver.}
\paragraph{Faraday flow}
Figure~\ref{fig:precision:faraday} shows error against runtime for a particle with $S=0.1$ and $R=\nicefrac{7}{9}$ in the Faraday flow, again for trajectories computed over \SI{1}{\second}.
All methods are more accurate than for the Bicklet jet, likely because the flow velocities in the Faraday flow are smaller.
For zero initial relative velocity, FD2 + IMEX2, FD4 + IMEX4 and Daitche~\reva{\cite{Daitche2013}} all perform very similarly with a slight advantage for FD4 + IMEX4 for higher accuracies.
For non-zero initial relative velocity, Daitche~\reva{\cite{Daitche2013}} again overtakes the FD methods with respect to speed of computation except for errors above $10^{-4}$, where FD2 becomes more efficient.
As for the Bickley jet, while Prasath et al.'s~\reva{\cite{PrasathEtAl2019}} method is highly accurate, it is too computationally costly to provide results in real-time.
\reva{Again, comparing the second- and fourth-order fully implicit methods with their implicit-explicit counterparts illustrates that the latter provides substantial gains in efficiency.}

\section{Conclusions}\label{section:conclusions}
The Maxey-Riley\reva{-Gatignol} equations (MRGE) describe the movement of inertial particles in a fluid and are used in a wide range of applications.
They contain an integral term,  called Basset history term, that makes them difficult to solve numerically.
Hence this term is often neglected, even though it has been shown in theory and experiments that this can lead to significant errors.
While direct numerical solvers for the full MRGE exist~\cite{Daitche2013}, they need to store all computed time steps to evaluate the integral term which can lead to memory problems, in particular when trajectories of many particles are computed concurrently~\cite{KopperEtAl2023}.
Approximations to the history term based on windowing are used but require careful calibration of parameters that do not necessarily generalize to very different flow fields~\cite{DorganAndLoth2007}.
A recent result by Prasath et al.~\cite{PrasathEtAl2019} shows that solving the full MRGE is equivalent to solving a one-dimensional diffusive partial differential equation with nonlinear boundary condition on a semi-infinite domain.
This allows to apply well established techniques for the numerical solution of partial differential equations to the MRGE and to avoid the memory demands of direct solvers.

However, solving the reformulated MRGE requires to deal with the semi-infinite computational domain.
Prasath et al.~\reva{\cite{PrasathEtAl2019}} propose a numerical approach based on polynomial expansions.
While accurate, this method is computationally expensive.
Instead, we adopt techniques developed by Koleva~\cite{Koleva2005} and Fazio and Jannelli~\cite{FazioAndJanelli2014} to develop finite difference discretizations of order two and four for the reformulated MRGE, resulting in a semi-discrete initial value problem.
Applying an explicit method would be inefficient due to the stiffness of the discrete Laplacian, whilst a fully implicit method would require nonlinear Newton solves in each step because of the nonlinear boundary condition.
Therefore, we propose to use an implicit-explicit RKM that treats the nonlinear boundary condition explicitly but the discrete Laplacian implicity.
This way, only linear solvers are required but the strict stability limit from the discrete Laplacian is avoided.

We test all methods, Prasath et al.'s polynomial expansion, Daitche's direct integrator and our FD2 and FD4 methods, for five different flow fields.
For three of them, a quiescient flow, a steady, inhomogeneous vortex and a non-steady, homogeneous oscillating background, analytical solutions are available.
For two unsteady, inhomogeneous fields, the Bickley jet and a Faraday flow interpolated from experimental data, we compute a reference  solution with Prasath et al.'s polynomial expansion approach with very high resolution.
Our numerical experiments suggest that:
\begin{enumerate}
 \itemsep0em
 \item Daitche's method is efficient for both zero and non-zero relative velocity and a good overall choice. 
 However, it can become unstable for very small density ratios and Stokes numbers. 
 \item Prasath's polynomial expansion based method is very accurate for both zero and non-zero relative velocity but computationally expensive.
 \item The finite difference method of order two is efficient for both zero and non-zero initial relative velocity but, particularly in the latter case, Daitche's method is significantly more efficient.
 \item The finite difference method of order four is efficient for zero initial relative velocity and could outperform both Daitche and FD2 + IMEX2 for the Faraday flow.
\end{enumerate}
Daitche and FD2 + IMEX2 could both simulate \SI{1}{\second} trajectories in significantly less than one second runtime, illustrating their potential to provide solutions to the MRGE in real-time.
For the Faraday flow, for example, FD2 + IMEX2 could compute a \SI{1}{\second}-trajectory with a relative discretization error of $7 \times 10^{-5}$ in \SI{8}{\milli\second}.
An optimized implementation in a compiled language plus effective use of vectorization and multi-threading would make simulations even faster.
This could potentially allow to use the full MRGE as a model in either closed-loop control or for model-based filtering to enable real-time particle tracking, for example of Lagrangian sensors that are too large to be considered perfect tracers~\cite{Bisgaard2020}.

\section*{Acknowledgments}
We are grateful to Prof. Alexandra von Kameke for providing us with the data from the experimental measurements of the Faraday flow and to Prof. Kathrin Padberg-Gehle for providing the Bickley jet example.

\appendix

\section{Boundary conditions for the 4th order upwind scheme}
\label{appendix:a}
Here we explain how to obtain a finite difference approximation of order four for the spatial derivative at the boundary $\bm{q}_x(0,t)$ in~\eqref{eq:MRSys}. 
The procedure is inspired by~\cite{malele2022}, where upwinding schemes for first and second derivatives for Robin boundary conditions are discussed.
Their general form is
\begin{align}
    a_0 (\partial_\xi q)_0 + a_1 (\partial_\xi q)_1 = \hat{b}_0 \left(\alpha q_0 - \gamma (\partial_x q)_0\right) + \frac{1}{d}\left\{ b_0q_0 + b_1q_1 + b_2q_2 + b_3q_3 \right\},
\end{align}
with $a_0,a_1,b_0,\hat{b}_0,b_1,\hat{b}_1,b_2,b_3\in\mathbb{R}$.
We use the relation between mesh derivatives~\eqref{Chain_rule:rel} to change the spatial derivative in the boundary condition from pseudo-space to the equidistant grid and obtain
\begin{align} \label{upwinding:scheme:1}
    a_0 (\partial_\xi q)_0 + a_1 (\partial_\xi q)_1 = \hat{b}_0 \left(\alpha q_0 - \frac{\gamma}{c} (\partial_\xi q)_0\right) + \frac{1}{d}\left\{ b_0q_0 + b_1q_1 + b_2q_2 + b_3q_3 \right\}.
\end{align}
We then expand all the elements of equation~\eqref{upwinding:scheme:1} around $\xi_0$ and find
\begin{align*}
    (\partial_\xi q)_0 :  & 0 &+& (\partial_\xi q)_0,\\
    (\partial_\xi q)_1 :  & 0 &+& (\partial_\xi q)_0 &+& d \; (\partial^2_\xi q)_0 &+& \frac{d^2}{2} (\partial^3_\xi q)_0 &+& \frac{d^3}{6} (\partial^4_\xi q)_0 &+& \frac{d^4}{24} (\partial^5_\xi q)_0,\\
    q_0 :  & q_0, \\
    q_1 :  & q_0 &+& d (\partial_\xi q)_0 &+& \frac{d^2}{2} (\partial^2_\xi q)_0 &+& \frac{d^3}{6} (\partial^3_\xi q)_0 &+& \frac{d^4}{24} (\partial^4_\xi q)_0 &+& \frac{d^5}{120} (\partial^5_\xi q)_0, \\
    q_2 :  & q_0 &+& 2d (\partial_\xi q)_0 &+& \frac{4d^2}{2} (\partial^2_\xi q)_0 &+& \frac{8d^3}{6} (\partial^3_\xi q)_0 &+& \frac{16d^4}{24} (\partial^4_\xi q)_0 &+& \frac{32d^5}{120} (\partial^5_\xi q)_0, \\
    q_3 :  & q_0 &+& 3d (\partial_\xi q)_0 &+& \frac{9d^2}{2} (\partial^2_\xi q)_0 &+& \frac{27d^3}{6} (\partial^3_\xi q)_0 &+& \frac{81d^4}{24} (\partial^4_\xi q)_0 &+& \frac{243d^5}{120} (\partial^5_\xi q)_0.
\end{align*}

By comparing coefficients we obtain the fourth order method
\begin{align*}
    q_0 : \quad & 0 &+& 0 &+& \alpha \hat{b}_0 &+& \frac{1}{d}b_0 &+& \frac{1}{d}b_1 &+& \frac{1}{d}b_2 &+& \frac{1}{d}b_3 &= 0, \\
    (\partial_\xi q)_0 : \quad & a_0 &+& a_1 &+& \frac{\gamma}{c}\hat{b}_0 &+& 0 &-& b_1 &-& 2b_2 &-& 3b_3 &= 0, \\
    (\partial^2_\xi q)_0 : \quad & 0 &+& da_1 &+& 0 &+& 0 &-& \frac{d}{2}b_1 &-& \frac{4d}{2}b_2 &-& \frac{9d}{2}b_3 &= 0, \\
    (\partial^3_\xi q)_0 : \quad & 0 &+& \frac{d^2}{2}a_1 &+& 0 &+& 0 &-& \frac{d^2}{6} b_1 &-& \frac{8d^2}{6} b_2 &-& \frac{27d^2}{6} b_3 &= 0, \\
    (\partial^4_\xi q)_0 : \quad & 0 &+& \frac{d^3}{6} a_1 &+& 0 &+& 0 &-& \frac{d^3}{24} b_1 &-& \frac{16d^3}{24} b_2 &-& \frac{81d^3}{24} b_3 &= 0, \\
    (\partial^5_\xi q)_0 : \quad & 0 &+& \frac{d^4}{24}a_1 &+& 0 &+& 0 &-& \frac{d^4}{120}b_1 &-& \frac{32d^4}{120}b_2 &-& \frac{243d^4}{120}b_3.
\end{align*}

We set $a_0 = a_1 = 1$.
This results in the coefficients
\begin{align*}
    \hat{b}_0 = -\frac{2c}{3\gamma} \quad b_0 = \frac{12\alpha dc-17\gamma}{18\gamma} \quad b_1 = \frac{1}{2} \quad b_2 = \frac{1}{2} \quad b_3 = -\frac{1}{18},
\end{align*}
and a leading error term of $\frac{d^4}{60}(\partial^5_\xi q)_0$.
In summary, the fourth order approximation at the boundary becomes
\begin{align*}
    (\partial_\xi q)_0 + (\partial_\xi q)_1 = -\frac{2c}{3\gamma} \left( \alpha q_0 - \frac{\gamma}{c} (\partial_\xi q)_0) \right) + \frac{1}{d} \left( \frac{12\alpha dc-17\gamma}{18\gamma}\right)q_0 + \frac{1}{2d}q_1 + \frac{1}{2d}q_2 - \frac{1}{18d}q_3,
\end{align*}
or written as a derivative with respect to $x$,
\begin{align*}
    c(\partial_x q)_0 + \frac{Nc}{N-1}(\partial_x q)_1 =& -\frac{2c}{3\gamma} \left( \alpha q_0 - \gamma (\partial_x q)_0 \right) \\
    &+\frac{1}{d} \left( \frac{12\alpha dc-17\gamma}{18\gamma} \right)q_0 + \frac{1}{2d}q_1 + \frac{1}{2d}q_2 - \frac{1}{18d}q_3.
\end{align*}
Finally, we use~\eqref{MRSys3} to get
\begin{align*}
    c(\partial_x q)_0 + \frac{Nc}{N-1}(\partial_x q)_1 =& -\frac{2c}{3\gamma} \left(f(q_0,x_0,t) - (\partial_t q)_0 \right) \\
    &+ \frac{1}{d} \left( \frac{12\alpha dc-17\gamma}{18\gamma} \right)q_0 + \frac{1}{2d}q_1 + \frac{1}{2d}q_2 - \frac{1}{18d}q_3.
\end{align*}

\section{Matrices for the fourth order finite difference method}
\label{appendix:b}

\begin{align*}
    M_1 :=
    \begin{bmatrix}
        c               & 0 & \frac{Nc}{N-1}           & 0 & 0               & 0              & 0 & \dots        & 0 \\
        0 & c               & 0 & \frac{Nc}{N-1}           & 0 & \ddots               & \ddots              & \ddots        & 0 \\
        c               & 0 & \frac{4Nc}{N-1}           & 0 & \frac{Nc}{N-2} & 0 & \ddots        & \ddots        & 0 \\
        0 & c               & 0 & \frac{4Nc}{N-1}           & 0 & \frac{Nc}{N-2} & 0        & \ddots        & 0 \\
        \vdots               & 0 & \frac{Nc}{N-1} & 0 & \frac{4Nc}{N-2}            & 0 & \frac{Nc}{N-3}  & \ddots        & 0 \\
        \vdots          & \ddots          & \ddots  & \ddots & \ddots & \ddots  & \ddots  & \ddots        & \vdots \\
        \vdots          & \ddots & \ddots & \frac{Nc}{4}  & 0 & \frac{4Nc}{3}       & 0 & \frac{Nc}{2} & 0 \\
        \vdots          & \ddots & \ddots & \ddots & \frac{Nc}{4}  & 0 & \frac{4Nc}{3}       & 0 & \frac{Nc}{2} \\
        0 & \ddots & \ddots & \ddots & \ddots & \frac{Nc}{3} & 0 & \frac{4Nc}{2}& 0 \\
        0 & 0 & \dots & \dots & \dots & 0 & \frac{Nc}{3} & 0 & \frac{4Nc}{2}
    \end{bmatrix},
\end{align*}

\begin{align*}
    M_2 :=
    \begin{bmatrix}
        c & 0 & 3\frac{Nc}{N-1} & 0 & 0 & 0 & 0 & 0 & \dots        & 0 \\
        0 & c & 0 & 3\frac{Nc}{N-1} & 0 & \ddots & \ddots & \ddots & \ddots        & 0 \\
        c & 0 & \frac{4Nc}{N-1} & 0 & \frac{Nc}{N-2} & 0 & \ddots & \ddots & \ddots & 0 \\
        0 & c & 0 & \frac{4Nc}{N-1} & 0 & \frac{Nc}{N-2} & 0 & \ddots & \ddots & 0 \\
        0 & 0  & \frac{Nc}{N-1} & 0 & \frac{4Nc}{N-2} & 0 & \frac{Nc}{N-3} & 0 & \ddots & 0 \\
        \vdots & \ddots & \ddots & \ddots & \ddots & \ddots & \ddots & \ddots & \ddots & \vdots \\
        \vdots & \ddots & \ddots & \ddots & \frac{Nc}{4} & 0 & \frac{4Nc}{3} & 0 & \frac{Nc}{2} & 0 \\
        \vdots & \ddots & \ddots & \ddots & 0 & \frac{Nc}{4} & 0 & \frac{4Nc}{3} & 0 & \frac{Nc}{2} \\
        0 & \ddots & \ddots & \ddots & \ddots & 0 & \frac{Nc}{3} & 0 & \frac{4Nc}{2} & 0 \\
        0 & 0 & \dots & \dots & \dots & \dots & 0 & \frac{Nc}{3} & 0 & \frac{4Nc}{2}
    \end{bmatrix},
\end{align*}

\begin{align*}
    B_1 := \frac{1}{d}
    \begin{bmatrix}
        \frac{12\alpha dc - 17\gamma}{18\gamma} & 0 & \nicefrac{1}{2} & 0 & \nicefrac{1}{2} & 0 & -\nicefrac{1}{18} & 0 & \dots & 0 \\
        0 & \frac{12\alpha dc - 17\gamma}{18\gamma} & 0 & \nicefrac{1}{2} & 0 & \nicefrac{1}{2} & 0 & -\nicefrac{1}{18} & \ddots & 0 \\
        -3 & 0 & 0 & 0 & 3 & 0 & 0 & 0 & \dots & 0 \\
        0 & -3 & 0 & 0 & 0 & 3 & 0 & 0 & \ddots & 0 \\
        \vdots & \ddots & \ddots & \ddots & \ddots & \ddots & \ddots & \ddots & \ddots & \vdots \\
        \vdots & \ddots & \ddots & \ddots & \ddots & 0 & -3 & 0 & 0 & 0\\
        0 & \dots & \dots & \dots & \dots & \dots & 0 & -3 & 0 & 0
    \end{bmatrix},
\end{align*}

\begin{align*}
    B_2 := \frac{1}{d}
    \begin{bmatrix}
        -\nicefrac{17}{6} & 0 & \nicefrac{3}{2} & 0 & \nicefrac{3}{2} & 0 & -\nicefrac{1}{6} & 0 & \dots & 0 \\
        0 & -\nicefrac{17}{6} & 0 & \nicefrac{3}{2} & 0 & \nicefrac{3}{2} & 0 & -\nicefrac{1}{6} & \dots & 0 \\
        -3 & 0 & 0 & 0 & 3 & 0 & 0 & 0 & \dots & 0 \\
        0 & -3 & 0 & 0 & 0 & 3 & 0 & \ddots & \ddots & 0 \\
        0 & 0 & -3 & 0 & 0 & 0 & 3 & 0 & \ddots & \vdots \\
        \vdots & \ddots & \ddots & \ddots & \ddots & \ddots & \ddots & \ddots & \ddots & \vdots \\
        \vdots & \ddots & \ddots & \ddots & \ddots & 0 & -3 & 0 & 0 & 0\\
        0 & \dots & \dots & \dots & \dots & \dots & 0 & -3 & 0 & 0
    \end{bmatrix},
\end{align*}.

\begin{align*}
    K_1 := 
    \begin{bmatrix}
        -\frac{2c}{3\gamma} f^{(1)}(\bm{q}_0(t),\bm{y}_0(t),t) \\
        -\frac{2c}{3\gamma} f^{(2)}(\bm{q}_0(t),\bm{y}_0(t),t) \\
        0 \\
        0 \\
        \vdots \\
        0
    \end{bmatrix}
\quad \text{and} \quad
    V_1 := \frac{2c}{3\gamma}
    \begin{bmatrix}
        \frac{\partial q^{(1)}(\xi(x),t)}{\partial t} \biggr\rvert_0 \\
        \frac{\partial q^{(2)}(\xi(x),t)}{\partial t} \biggr\rvert_0 \\
        0 \\
        0 \\
        \vdots \\
        0
    \end{bmatrix}.
\end{align*}

% \section{Numerical convergence rates}
% \input{appC}

\bibliography{refs}

\end{document}